\documentclass[12pt,a4paper]{amsart}
\usepackage[utf8]{inputenc}
\usepackage{amsfonts}
\usepackage{amssymb}
\usepackage{bbm}
\usepackage{comment}
\usepackage{dsfont}
\usepackage{amsthm}
\usepackage{amsmath}
\usepackage{commath}
\usepackage{url}
\usepackage{epsfig}
\usepackage{adjustbox}
\usepackage[ruled,vlined,linesnumbered]{algorithm2e}
\usepackage{cool}
\usepackage[table]{xcolor}
\usepackage{tabularx}
\usepackage[a4paper]{geometry}
\usepackage{caption}
\usepackage{booktabs}
\usepackage{subcaption}
\usepackage{stmaryrd}
\usepackage{graphicx,accents}
\usepackage{enumerate}

\usepackage{a4wide}
\usepackage{todonotes}
\usepackage{scrextend}
\usepackage{mathtools}

\newcommand{\ls}{\leqslant}
\newcommand{\gs}{\geqslant}

\makeatletter
\DeclareRobustCommand{\cev}[1]{%
  \mathpalette\do@cev{#1}%
}
\newcommand{\do@cev}[2]{%
  \fix@cev{#1}{+}%
  \reflectbox{$\m@th#1\vec{\reflectbox{$\fix@cev{#1}{-}\m@th#1#2\fix@cev{#1}{+}$}}$}%
  \fix@cev{#1}{-}%
}
\newcommand{\fix@cev}[2]{%
  \ifx#1\displaystyle
    \mkern#23mu
  \else
    \ifx#1\textstyle
      \mkern#23mu
    \else
      \ifx#1\scriptstyle
        \mkern#22mu
      \else
        \mkern#22mu
      \fi
    \fi
  \fi
}

\newtheorem{thm}{Theorem}
\newtheorem{cor}[thm]{Corollary}
\newtheorem{lem}[thm]{Lemma}
\newtheorem{prop}[thm]{Proposition}
\newtheorem{rmk}{Remark}

\title[On the asymptotic behaviour of Sudler products]{On the asymptotic behaviour of Sudler products for badly approximable numbers}
\author{Manuel Hauke}

\begin{document}

\maketitle
\begin{abstract}
    Given a badly approximable number $\alpha$, we study the asymptotic
    behaviour of the Sudler product defined by 
    $P_N(\alpha) = \prod_{r=1}^N 2 \lvert \sin \pi r \alpha \rvert$.
    We show that $\liminf_{N \to \infty} P_N(\alpha) = 0$ and $\limsup_{N \to \infty} P_N(\alpha)/N = \infty$ whenever the
    sequence of partial quotients in the continued fraction expansion of $\alpha$ exceeds $7$ infinitely often. This improves results obtained by Lubinsky for the general case, and by Grepstad, Neumüller and Zafeiropoulos for the special case of quadratic irrationals. Furthermore, we prove that this threshold value $7$ is optimal,
    even when restricting $\alpha$ to be a quadratic irrational, which gives a negative answer to a question of the latter authors.
\end{abstract}

\section{Introduction and statement of results}

\subsection{Introduction}

For $\alpha \in \mathbb{R}$ and $N$ a natural number, the Sudler product is defined as
\begin{equation*}
P_N(\alpha) := \prod_{r=1}^{N} 2 \left\lvert \sin \pi r\alpha \right\rvert.
\end{equation*}
Note that by $1$-periodicity of $P_N(\alpha)$ and the fact that $P_N(\alpha) = 0$ for rational $\alpha$ and $N$ sufficiently large, the asymptotic analysis of such products can be restricted
to irrational numbers $\alpha \in [0,1]$.\\
Sudler products appear in many different areas of mathematics that include, among others, restricted partition functions \cite{sudler}, KAM theory \cite{kam} and Padé approximants \cite{lubinsky_pade}, and were used in the context of almost Mathieu operators in the solution of the Ten Martini Problem by Avila and Jitomirskaya \cite{aj}. 
Recently, Aistleitner and Borda \cite{quantum_invariants,zag_conj} established a connection between Sudler products and the work of Bettin and Drappeau \cite{bd1,bd2,bd3} on the order of magnitude of the Kashaev invariant of certain hyperbolic knots, following the work of Zagier \cite{zag}.\par

Writing $\| P_N \|_{\infty} = \max_{0 < \alpha < 1} P_N(\alpha)$, Erd\H os and Szekeres \cite{erdos_szekeres} claimed that the limit
$\lim_{N \to \infty} \| P_N \|_{\infty}^{1/N}$ exists and equals a value between $1$ and $2$, without formally proving it. This was done by Sudler \cite{sudler} and Wright \cite{wright} who showed that
$\lim_{N \to \infty} \| P_N \|_{\infty}^{1/N} = C \approx 1.22$.
Inspired by this, the order of growth of Sudler products was extensively examined from a metric point of view. For more results in this area, we refer
the reader to \cite{atkinson,bbr,bell,bglly,bc,fh,kol,kol2}.\\

In this paper, we consider the pointwise behaviour of Sudler products.
First studied by Erd\H os and Szekeres \cite{erdos_szekeres}, 
it was proven that

\begin{equation}\label{asym=0}\liminf_{N \to \infty} P_N(\alpha) = 0, \quad \limsup_{N \to \infty} P_N(\alpha) = \infty\end{equation}
holds for almost every $\alpha$.
In contrast to the result on $\lim_{N \to \infty} \| P_N \|_{\infty}^{1/N}$, Lubinsky and Saff \cite{ls} showed that for almost every $\alpha$, \mbox{$\lim_{N \to \infty} P_N(\alpha)^{1/N} = 1$}.
For distributional results for varying $N$ and fixed $\alpha$, see e.g. \cite{borda,hauke_density}.

Mestel and Verschueren \cite{mestel} examined the behaviour of $P_N(\phi)$ where $\phi = \frac{\sqrt{5}-1}{2} = [0;1,1,1,\ldots]$ is the fractional part of the Golden Ratio.
 They showed that the limit along the Fibonacci sequence $\lim\limits_{n \to \infty} P_{F_n}(\phi)$ exists.
 As the Fibonacci numbers are the denominators of the continued fraction convergents of $\phi$, this hints at a connection between the Sudler product of $\alpha$ and its Diophantine approximation properties. Indeed, this was established in \cite{quantum_invariants,tech_zaf,grepstadII}, especially for the case where $\alpha$
 is a quadratic irrational.
 \par

Returning to \eqref{asym=0}, Lubinsky \cite{lubinsky} showed that 
$\liminf_{N \to \infty} P_N(\alpha) = 0$ holds for any $\alpha \in \mathbb{R}$ that has unbounded 
continued fraction coefficients.
In fact, he proved that there is a finite cutoff value $K$ (in 
\cite{grepstad_survey} it was shown that the proof of Lubinsky actually gives a value of $K \approx e^{800}$) such that whenever an irrational $\alpha = [a_0;a_1,a_2,\ldots]$ fulfills $\limsup_{n \to \infty} a_n \gs K$, then $\liminf_{N \to \infty} P_N(\alpha) = 0$.
Lubinsky conjectured that 
$\liminf_{N \to \infty} P_N(\alpha) = 0$ holds for all $\alpha$,
which was disproven by Grepstad, Kaltenböck and Neumüller \cite{grepstad_neum}, as they 
proved that $\liminf_{N \to \infty} P_N(\phi) > 0.$\par

Aistleitner, Technau and Zafeiropoulos \cite{tech_zaf} found a close connection between the behaviour of $\liminf\limits_{N \to \infty} P_N(\alpha)$ and $\limsup\limits_{N \to \infty} \frac{P_N(\alpha)}{N}$, which was generalized by Aistleitner and Borda in  \cite{quantum_invariants}: They showed that for badly approximable numbers\footnote{Although \eqref{liminf_limsup_equiv} is only stated for quadratic irrationals, 
the result follows for any badly approximable number as Aistleitner and Borda show that for $0 \ls N \ls q_k$, we have
$\log P_N(\alpha) + \log P_{q_k-N-1}(\alpha) = \log q_k + \mathcal{O}(\log (1 + \max\limits_{k \in \mathbb{N}} a_k)).$
} $\alpha$, we have
\begin{equation}\label{liminf_limsup_equiv}\liminf_{N \to \infty} P_N(\alpha) = 0 \; \Longleftrightarrow \; \limsup_{N \to \infty} \frac{P_N(\alpha)}{N} = \infty.\end{equation}

It was shown in \cite{tech_zaf} that for $\beta(b) := [0;\overline{b}]$, we have
$\liminf\limits_{N \to \infty} P_N(\beta) > 0$ and $\limsup\limits_{N \to \infty}      \frac{P_N(\beta)}{N} < \infty$, if and only if $b \ls 5$. Recently, the sufficient condition for $\liminf\limits_{N \to \infty} P_N(\alpha) = 0$ and $\limsup\limits_{N \to \infty} \frac{P_N(\alpha)}{N} = \infty$ was generalized by  Grepstad, Neumüller and Zafeiropoulus \cite{grepstadII} to arbitrary quadratic irrationals $\alpha$. They showed that if $\alpha = [a_0;a_1,\ldots,a_{p},\overline{a_{p+1},\ldots,a_{p + \ell}}]$ and
$a_K := \limsup_{n \to \infty} a_n \gs 23$, then $\liminf\limits_{N \to \infty} P_N(\alpha) = 0$ and $\limsup\limits_{N \to \infty} \frac{P_N(\alpha)}{N} = \infty$, also showing that this value can be decreased to $22$ when the (shortest) period length $\ell$ is even.
In the following theorem, we prove that this value can be decreased to $7$, even when considering
arbitrary badly approximable numbers.

\begin{thm}\label{K_l<=7}
Let $\alpha$ be an irrational number with continued fraction expansion $\alpha = [a_0;a_1,a_2,\ldots]$ that satisfies
\[\limsup_{n \to \infty} a_n \gs 7.\]
Then 
\begin{equation*}
\liminf_{N \to \infty} P_N(\alpha) = 0, \quad \limsup_{N \to \infty} \frac{P_N(\alpha)}{N} = \infty.
\end{equation*}
\end{thm}

 Based on numerical evidence, Grepstad, Neumüller and Zafeiropoulus \cite{grepstadII} speculated that $\liminf_{N \to \infty}P_N(\alpha) = 0$ holds for any quadratic irrational that fulfills $\limsup_{n \to \infty} a_n(\alpha) \gs 6$, which would be in accordance with the result obtained by Aistleitner, Technau and Zafeiropoulus \cite{tech_zaf} for the special case of period length $\ell = 1$. We show that this remains true when $\ell =2$.

\begin{thm}\label{l=2}
Let $\alpha$ be a quadratic irrational with continued fraction expansion\\
\mbox{$\alpha = [a_0;a_1,\ldots,a_p,\overline{a_{p+1},a_{p+2}}]$}
and assume that $a_K = \max\{a_{p+1},a_{p+2}\} \gs 6$.
Then we have

\begin{equation*}
\liminf_{N \to \infty} P_N(\alpha) = 0, \quad \limsup_{N \to \infty} \frac{P_N(\alpha)}{N} = \infty.
\end{equation*}
\end{thm}

However, the following theorem proves that in general, the condition $\limsup_{n \to \infty}a_n(\alpha) \gs 6$ does not imply $\liminf_{N \to \infty} P_N(\alpha) = 0$, even when restricting to quadratic irrationals.

\begin{thm}\label{counterex_l3}
Let $\alpha = [0;\overline{6,5,5}]$ or $\alpha = [0;\overline{5,4}]$. Then we have

\begin{equation*}
\liminf_{N \to \infty} P_N(\alpha) > 0, \quad \limsup_{N \to \infty} \frac{P_N(\alpha)}{N} < \infty.
\end{equation*}
\end{thm}

Clearly, this gives a negative answer to the question raised by Grepstad, Neumüller and Zafeiropoulus. Further, note that the combination of Theorems \ref{K_l<=7} (respectively Theorem \ref{l=2}) and \ref{counterex_l3} leads to the result that the threshold value for $\limsup_{n \to \infty} a_n(\alpha)$ in Theorem \ref{K_l<=7} (respectively Theorem \ref{l=2}) is indeed optimal. This strongly improves the current best bounds which were $23$ for quadratic irrationals and $K \approx e^{800}$ for arbitrary irrationals.

\subsection{Directions for further research}
\begin{itemize}
    \item
     Fixing a period length $\ell \gs 1$, what is the smallest value for $a_K = \limsup_{n \to \infty}a_n$ such that any $\alpha = [a_0;a_1,\ldots,a_{p},\overline{a_{p+1},\ldots,a_{p + \ell}}]$ with $\ell$ chosen minimal,
     fulfills \eqref{asym=0}? Denoting this optimal value by
     $K_{\ell}$, the results from \cite{tech_zaf} and Theorems \ref{K_l<=7}\,--\,\ref{counterex_l3} show that $K_1 = 6, K_2 = 6, K_3 = 7, K_{\ell} \ls 7$ for $\ell \gs 4$. We believe that $K_{\ell} = 7$ for any $\ell \gs 4$, and that this can be proven by considering the irrationals
     $\alpha_{\ell} = [0;\overline{6,\underbrace{5,\ldots,5}_{\ell-1 \text{ times}}}]$.\\
     In principle, this should be possible to prove with the methods applied in Theorem \ref{counterex_l3} for any fixed $\ell$, but both the combinatorial and computational effort increases for larger $\ell$. Also for $\ell \to \infty$, one needs to establish an additional argument.
     \vspace{3mm}
     
     \item We are also interested in the dual problem in the following sense: what is the maximal value of $a_K$ (or which other property) that guarantees $\liminf_{N \to \infty} P_N(\alpha) > 0, \quad \limsup_{N \to \infty} \frac{P_N(\alpha)}{N} < \infty$? For quadratic irrationals with $\ell = 1$, \cite{tech_zaf} proves that this is the case for $a_K = 5$. However, for general quadratic irrationals (and thus badly approximable integers), the optimal value cannot exceed $3$, since we show in a byproduct of the proof of Theorem \ref{l=2} that $\alpha = [0;\overline{4,1}]$ still fulfills $\liminf_{N \to \infty} P_N(\alpha) = 0$. We believe that the value $a_K = 3$ is optimal, however, our tools are not strong enough to prove this: the criterion used to prove Theorem \ref{K_l<=7} uses just a sufficient, but not necessary condition for $\alpha$ to fulfill $\liminf_{N \to \infty} P_N(\alpha) = 0$, and the methods in Theorems \ref{l=2} and \ref{counterex_l3} can only be applied to a fixed irrational $\alpha$ (modulo preimages under Gauss map iterations), and not to infinitely many $\alpha$ at the same time. Thus again, we could only prove such results for quadratic irrationals with fixed period length $\ell$.
     
     \item Lubinsky \cite{lubinsky} showed that for any badly approximable $\alpha$, there exist 
     $c_1(\alpha), c_2(\alpha)$ such that
     \begin{equation}\label{c1c2}N^{-c1} \ll P_N(\alpha) \ll N^{c_2}.\end{equation}
     Note that this is in contrast to the almost sure behaviour (with respect to Lebesgue measure). For almost every $\alpha$ and every $\varepsilon > 0$, we have
     \[\log N \log \log N \ll \log P_N(\alpha) \ll \log N (\log \log N)^{1+ \varepsilon},\]
     which can be deduced from \cite{borda} and \cite{lubinsky}.
     It was shown in \cite{quantum_invariants} that
     if $C_1(\alpha),C_2(\alpha)$ denote the infimum of all $c_1,c_2$ fulfilling \eqref{c1c2}, we have $C_2(\alpha) = C_1(\alpha) +1$. For more results in this direction, we refer the reader to \cite{quantum_invariants,other_aist_borda} where several statements about estimates on $C_1(\alpha), C_2(\alpha)$ are proven. 
    If $\liminf_{N \to \infty} P_N(\alpha) > 0$, then clearly we have $C_1(\alpha) = 0$. 
    To the best of our knowledge, the precise value of $C_1(\alpha)$ is not known for 
    any $\alpha$ that fulfills $\liminf_{N \to \infty} P_N(\alpha) = 0$.
     \vspace{3mm}
     
     \item In \cite{hauke_extreme}, we proved that for $\beta = [0;\overline{b}]$, $b \ls 5$, and $q_{n-1} \ls N \ls q_n$, we have that
     $P_{q_{n -1}}(\beta) \ls P_N(\beta) \ls P_{q_n-1}(\beta)$. We are interested in whether a similar behaviour holds for any badly approximable number $\alpha$ (with the previous inequality adjusted to suitable subsequences $q_{n_k}$), provided that $n$ is chosen sufficiently large, and $\alpha$ fulfills $\liminf_{N \to \infty} P_N(\alpha) > 0$.
\end{itemize}

\subsection{Structure of the paper}
The rest of the article is structured as follows: In Section \ref{prereq},
we recall basic properties of continued fractions that are used in this article
and review already established results on perturbed Sudler products and, if $\alpha$ is a quadratic irrational, their limiting behaviour along subsequences.
In Section \ref{mainproof}, we introduce a sufficient condition for badly approximable
$\alpha$ to deduce $\liminf_{N \to \infty} P_N(\alpha) = 0$, which can be applied for
infinitely many irrationals at the same time, and using this tool, we prove Theorem \ref{K_l<=7}.
In Sections \ref{secl2} and \ref{counterex_sec}, we prove Theorems \ref{l=2} respectively \ref{counterex_l3} by explicitly computing values of finitely many limit functions evaluated
at certain perturbations. This article is structured in a way such that the proofs in Sections \ref{mainproof}\,--\,\ref{counterex_sec} are almost independent of each other, only relying on statements proven in Section \ref{prereq}.

\section{Prerequisites}\label{prereq}

\subsection{Notation}

Given two functions $f,g:(0,\infty)\to \mathbb{R},$ we write $f(t) = \mathcal{O}(g(t)), f \ll g$ or $g \gg f$ when
$\limsup_{t\to\infty} \frac{|f(t)|}{|g(t)|} < \infty$. Any  dependence of the value of the limsup above on potential parameters is denoted by the appropriate subscripts.
Given a real number $x\in \mathbb{R},$ we write $\{x\}$ for the fractional part of $x$ and $\|x\|=\min\{|x-k|: k\in\mathbb{Z}\}$ for the distance of $x$ from its nearest integer.
We denote the characteristic function of a relation $R$ by $\mathds{1}_R$ and understand
the value of empty sums and products as $0$ respectively $1$.

\subsection{Continued fractions}

For convenience of the reader, we recall here some well-known facts on continued fractions that are used in this paper. For a more detailed background, see the classical literature e.g. \cite{all_shall,rock_sz,schmidt}. Furthermore, we define notations in context of continued fractions that are particularly useful here.
Every irrational $\alpha$ has a unique infinite continued fraction expansion $[a_0;a_1,...]$ with  convergents $p_k/q_k = [a_0;a_1,...,a_k]$ fulfilling the recursions 
\begin{equation*}
p_{k+1} = p_{k+1}(\alpha) = a_{k+1}p_k + p_{k-1}, \quad q_{k+1} = q_{k+1}(\alpha) = a_{k+1}q_k + q_{k-1}, \quad k \gs 1\end{equation*}
with initial values $p_0 = a_0,\; p_1 = a_1a_0 +1,\; q_0 = 1,\; q_1 = a_1$.
One can deduce from these recursions that $q_k$ grows exponentially fast in $k$; in particular,
we have for any $k,j \in \mathbb{N}$
\begin{equation}\label{exp_grow}
    \frac{q_{k}}{q_{k+j}} \gg_{\alpha} C^j
\end{equation}
where $0 < C(\alpha) < 1$.
Defining 
\[\delta_k := \lVert q_k \alpha\rVert,\quad 
\cev{\alpha}_{k} = [0;a_{k},a_{k-1},\ldots,a_1], \quad
\vec{\alpha}_{k} = [0;a_k,a_{k+1},\ldots],\quad k \gs 1,
\]
 we have 

\begin{align}\label{alternating_approx}
    q_k\alpha 
&\equiv (-1)^k\delta_k \pmod 1,
\\q_k\delta_k &= \frac{1}{a_{k+1} + \vec{\alpha}_{k+2} + \cev{\alpha}_{k}}.\label{qk_delta_k}
\end{align}

We call an irrational number $\alpha$ \textit{badly approximable} if $\limsup_{n \to \infty} a_n(\alpha) < \infty$. For fixed $a \in \mathbb{N}$, we define

\[\mathcal{B}(a,M) := \{\alpha  = [0;a_1,a_2,\ldots] \in \{0,1\}:
 \max_{i \gs M+1} a_i = a \}.
\]

If $\beta = [a_0;a_1,a_2,\ldots]$ is badly approximable with $\limsup_{n \to \infty} a_n(\beta) = a$, then there exists a $K_0 \in \mathbb{N}$ such that
$T^{K_0}(\beta) \in \mathcal{B}(a,M)$, where $T$ denotes the Gauss map defined by \[T(x) = \begin{cases}0 &\text{ if } x= 0\\ \left\{\frac{1}{x}\right\} \mod 1&\text{ if } 0 < x \ls 1.
\end{cases}\]

Fixing an irrational $\alpha = [a_0;a_1,...]$, the Ostrowski expansion of a non-negative integer $N$ is the unique representation

\begin{equation*}
N = \sum_{\ell = 0}^k b_{\ell}q_{\ell} \quad \text{ where }
b_{k+1} \neq 0,\;\; 0 \ls b_0 < a_1, \;\;  0 \ls b_{\ell} \ls a_{\ell+1} \text { for } \ell \gs 1,
\end{equation*}
with the additional rule that
$b_{\ell-1} = 0$ whenever $b_{\ell} = a_{\ell+1}$.
If $\alpha = \frac{p}{q}$ is a rational number (with $p,q$ coprime), then we have
$\alpha = [a_0;a_1,\ldots,a_{k}]$ for some $k \in \mathbb{N}$ with $a_k > 1$. For
$N < q = q_k$, the Ostrowski expansion is defined as in the irrational setting.\vspace{5mm}\par

We turn our attention to the case where $\alpha$ is a quadratic irrational,
which will be of particular interest in Theorems \ref{l=2} and \ref{counterex_l3}.
We define
\[ \mathcal{Q}(a) := \{\alpha = [0;\overline{a_1,\ldots,a_{\ell}}] \in \{0,1\}: \ell \in \mathbb{N}, \max_{1 \ls i \ls \ell} a_i = a\}\]
where $(\overline{a_1,\ldots,a_{\ell}}) := 
(a_1,\ldots,a_{\ell},a_1,\ldots,a_{\ell},a_1,\ldots,a_{\ell},\ldots)$.
As before, each quadratic irrational $\beta = [a_0;a_1,\ldots,a_p,\overline{a_{p+1},\ldots, a_{p+\ell}}]$ with $\max_{p+1 \ls k \ls \ell + p} a_k = a$ fulfills $T^p(\beta) \in \mathcal{Q}(a)$. We will see later that the arguments we use are invariant under the Gauss map, so it suffices to examine
quadratic irrationals $\alpha =  [0;\overline{a_1,\ldots,a_{\ell}}] \in \mathcal{Q}(a)$. For fixed $\ell$, we denote by $[k]$ the smallest non-negative residue of $k \mod \ell$.
 Furthermore, we adopt the following notation from \cite{grepstadII}. For ${\bf a} = (a_1,\ldots,a_{\ell})$, we define the permutation operators $\tau_r, \sigma_r$
for $0 \ls r \ls \ell-1$ by

\[
\tau_r({\bf a}) := (a_{r+1},\ldots,a_{\ell},a_1,\ldots,a_r), 
\]
respectively
\[
\sigma_r({\bf a}) := (a_{r-1},\ldots,a_1,a_{\ell},\ldots,a_r) \text{ if } r \gs 2,\]
and $\sigma_0({\bf a}) := (a_{\ell-1},\ldots,a_1,a_{\ell}),\sigma_1({\bf a}) := (a_{\ell},\ldots,a_1).
$
Its corresponding quadratic irrationals will be denoted by
\begin{equation}\label{alphatau}\alpha_{\tau_r} := [0;\overline{a_{r+1},\ldots,a_{\ell},a_1,\ldots,a_r}], \quad
\alpha_{\sigma_r} := [0;\overline{a_{r-1},\ldots,a_1,a_{\ell},\ldots,a_r}].
\end{equation}
With this notation, we obtain that for every $0 \ls r \ls \ell-1$
\begin{align}\label{quadr_limit1}
    \lim_{m \to \infty} \frac{q_{m\ell + r -1}}{q_{m\ell +r}} = \alpha_{\sigma_r}, \quad
    \lim_{m \to \infty} q_{m\ell + r}\delta_{m\ell + r} = \frac{1}{a_{r+1}+ \alpha_{\tau_{[r+2]}}+ \alpha_{\sigma_r}} =: C(r). 
\end{align}

\subsection{Perturbed Sudler products}
All results in this paper rely on a decomposition approach that was implicitly used already in \cite{grepstad_neum} to prove $\liminf\limits_{N \to \infty}P_N(\phi) > 0$, and more explicitly in later literature (see e.g. \cite{quantum_invariants,zag_conj,other_aist_borda,tech_zaf,grepstadII,hauke_density,hauke_extreme}). The Sudler product $P_N(\alpha)$ is
decomposed into a finite product with factors of the form

 \begin{equation}
    \label{shifted_sudler_beta}
P_{q_n}(\alpha,\varepsilon) := \prod_{r=1}^{q_n} 2 \left\lvert \sin\Big(\pi\Big(r\alpha + (-1)^{n}\frac{\varepsilon}{q_n}\Big)\Big)\right\rvert,
\end{equation}

a fact that is summarized by the following proposition:
\begin{prop}\label{prop_shifted}
Let 
$N = \sum_{i=0}^{n} b_{i}q_i(\alpha)$ be the Ostrowski expansion of an arbitrary integer $q_n \ls N < q_{n+1}$.
For $0 \ls i \ls n$ and $k \in \mathbb{N}$, we define
\begin{equation*}
\varepsilon_{i,k}(N)
:= q_i\left(k\delta_i + \sum_{j=1}^{n-i} (-1)^jb_{i+j} \delta_{i+j}\right)
\end{equation*}
and
\begin{equation}\label{def_Ki}K_i(N) := \prod_{i=0}^{n}\prod_{c_i= 1}^{b_i(N)-1} P_{q_i}(\alpha,\varepsilon_{i,c_i}(N))
\cdot P_{q_{i-1}}(\alpha,\varepsilon_{i-1,0}(N))^{\mathds{1}_{[b_{i-1}(N) \neq 0]}},\end{equation}
where $b_{-1} := 0$. Then we have

\begin{equation}\label{shifted_prod_1}
    P_N(\alpha) = \prod_{i=0}^{n}\prod_{c_i= 0}^{b_i-1} P_{q_i}(\alpha,\varepsilon_{i,c_i}(N)),
\end{equation}
as well as
\begin{equation}\label{shifted_prod_2}
    P_N(\alpha) = P_{q_n}(\alpha) \cdot \prod_{i = 1}^n K_i(N).
\end{equation}
\end{prop}

\begin{proof}
For $N = \sum_{i=0}^n b_{i}q_i(\alpha)$ fixed, we
define 
\begin{equation*}
M_{i,k} = M_{i,k}(N) = \sum_{j=i+1}^{n} b_jq_j + kq_i, \quad\quad i = 0,\ldots, n, \quad k = 0,\ldots,b_i-1.\end{equation*}
Using the approach of \cite{quantum_invariants,tech_zaf,grepstadII}, we obtain 

\begin{align*}
P_N(\alpha) = \prod_{i=0}^{n}\prod_{c_i= 0}^{b_i-1} P_{q_i}(\alpha,\tilde{\varepsilon}_{i,c_i})
\end{align*}
where by periodicity of the sine, any $\tilde{\varepsilon}_{i,c_i} = \tilde{\varepsilon}_{i,c_i}(N)$ that fulfills the relation

\begin{align}
\label{form_of_eps}\frac{(-1)^i\tilde{\varepsilon}_{i,c_i}}{q_i} &\equiv M_{i,c_i}\alpha \pmod 1
 \end{align}
 can be chosen. By \eqref{alternating_approx},
we have

\[M_{i,c_i}\alpha =  c_iq_i\alpha + \sum_{j=1}^{n-i} b_{i+j}q_{i+j}\alpha
\equiv (-1)^i\left(c_i\delta_i + \sum_{j=1}^{n-i} (-1)^jb_{i+j}\delta_{i+j}\right) \pmod 1,
\]
thus $\varepsilon_{i,c_i}(N)$ fulfills \eqref{form_of_eps}. This proves \eqref{shifted_prod_1}, whereas \eqref{shifted_prod_2} follows by a regrouping of the products.
\end{proof}

Mestel and Verschueren \cite{mestel} showed that 
$P_{F_n}(\phi) = P_{F_n}(\phi,0)$ converges to some positive constant $C_1$.
Aistleitner, Technau and Zafeiropoulus showed that for quadratic irrationals with period length $1$, the limit
$G(\alpha,\varepsilon) := \lim_{n \to \infty}P_{q_n}(\alpha,\varepsilon)$ exists and that the convergence is uniform on compact intervals.
The authors also give a rather long, closed expression for $G(\alpha,\varepsilon)$, that helped them to approximately compute $G(\alpha,\varepsilon)$ and in particular, $\lim_{N \to \infty} P_{q_n}(\alpha) = G(\alpha,0)$.
This result was generalized in \cite{grepstadII} to hold for quadratic irrationals with arbitrary period length in the following way:\\

{\bf Theorem A} (Grepstad, Neumüller, Zafeiropoulus \cite{grepstadII}). 
Let $\alpha := [0;\overline{a_1,\ldots,a_{\ell}}]$ and $P_{q_n}(\alpha,\varepsilon)$ be as in \eqref{shifted_sudler_beta}.
Then for each $1 \ls r \ls \ell$, $P_{q_{m\ell + k}}(\alpha,\varepsilon)$ converges locally uniformly to a function $G_{r}(\alpha,\varepsilon)$. The functions $G_r(\alpha,\cdot)$ are continuous and $C^{\infty}$ on every interval where they are non-zero.\vspace{4mm}

The limit functions play a crucial role in the sufficient condition used in \cite{grepstadII} to obtain their results. The following result is a straightforward generalization of \cite[Lemma 1]{tech_zaf} for quadratic irrationals with arbitrary period length $\ell$\vspace{4mm}:

{\bf Theorem B} (Grepstad, Neumüller, Zafeiropoulus \cite{grepstadII}). \label{G0}
Let $\beta = [0;a_1,\ldots a_p,\overline{a_{p+1},\ldots,a_{p +\ell}}]$ 
be an arbitrary quadratic irrational,
$\alpha = T^p(\beta) = [0;\overline{a_1,\ldots,a_{\ell}}]$ and $G_{r}(\alpha,\varepsilon)$ as in Theorem A. If there exists some $1 \ls r \ls \ell$ such that $G_{r}(\alpha,0) < 1$, then 
\begin{equation*}
\liminf_{N \to \infty} P_N(\beta) = 0, \quad \limsup_{N \to \infty} \frac{P_N(\beta)}{N} = \infty.
\end{equation*}
The following Lemma is helpful for any asymptotical analysis of Sudler products using the decomposition into perturbed products. It implies that the behaviour of finitely many shifted products is negligible for the asymptotic order of magnitude.

\begin{lem}\label{gg1}
Let $\alpha \in \mathcal{Q}(a), N \in \mathbb{N}$ and $\varepsilon_{i,c_i}(N)$ as before.
Then we have
\[1 \ll_{\alpha} P_{q_k}(\alpha,\varepsilon_{i,c_i}(N)) \ll_{\alpha} 1\]
uniformly in $N$ and $k$.
\end{lem}

\begin{proof}
The fact that $P_{q_k}(\alpha,\varepsilon_{i,c_i}(N)) \gg_{\alpha} 1$ is proven in \cite[Lemma 3]{quantum_invariants}, whereas $P_{q_k}(\alpha,\varepsilon_{i,c_i}(N)) \ll_{\alpha} 1$ follows immediately by the convergence to the uniformly bounded limit functions $G_r$.
\end{proof} 

\section{Proof of Theorem \ref{K_l<=7}}\label{mainproof}
As we want to prove a statement on arbitrary badly approximable numbers, the usage of limit functions from Theorem A is not possible.
However, the following lemma can be seen as a generalization of Theorem B 
for arbitrary badly approximable irrationals $\alpha$.

\begin{lem}\label{thmb_equiv}
Let $\alpha$ be a badly approximable number
and assume that there exist 
real numbers $\delta > 0, 0 < c < 1$ such that for infinitely many $k$ we have

\begin{equation}\label{shifted<1}\sup_{\varepsilon \in (-\delta,\delta)} P_{q_k}(\alpha,\varepsilon) < c.\end{equation}
Then 

\[\liminf_{N \to \infty} P_N(\alpha) = 0, \quad \limsup_{N \to \infty} \frac{P_N(\alpha)}{N} = \infty.\]
\end{lem}

\begin{proof}[Proof of Lemma \ref{thmb_equiv}]
The proof follows along the lines of \cite[Lemma 1]{tech_zaf} and \cite[Theorem 2]{grepstadII}. As there are infinitely many $k$ fulfilling \eqref{shifted<1}, 
$\delta_k \ll_{\alpha} \frac{1}{q_k}$
and $q_k$ grows exponentially in $k$,
we can extract an increasing subsequence $k_j$ such that the following holds for every $j \in \mathbb{N}$:

\begin{itemize}
    \item $\sup_{\varepsilon \in (-\delta,\delta)} P_{q_{k_j}}(\alpha,\varepsilon) < c < 1.$
    \item $q_{k_j} \gs 2q_{k_{j-1}}$, $j \gs 2$.
    \item $\delta_{k_j} < \frac{\delta}{4q_{k_{j-1}}}$.
\end{itemize}
Now we define $N_{j} := \sum\limits_{i =1}^j q_{k_j}, j \gs 1$. By Proposition \ref{prop_shifted},
we obtain 
\[P_{N_j}(\alpha) = \prod_{i=1}^{j}P_{q_{n_j}}(\alpha,\varepsilon_{n_j,0}(N)).
\]
Observe that
\[
\lvert \varepsilon_{n_j,0}(N) \rvert 
\ls q_{n_j}\left(
\sum_{\ell = 1}^{\infty} \delta_{n_{j+\ell}}
\right)
\ls \frac{\delta}{4}\left(\sum_{\ell = 0}^{\infty}\frac{1}{2^{\ell}}\right) = \frac{\delta}{2},
\]
so we obtain

\[P_{N_{k_j}}(\alpha) = 
\prod_{i=1}^j  P_{q_k}(\alpha,\varepsilon_{N_{k_j}}) < c^j,
\]
hence with $j \to \infty$, we obtain 
$\liminf_{N \to \infty} P_N(\alpha) = 0$. The statement follows now immediately from
\eqref{liminf_limsup_equiv}.
\end{proof}

In order to apply Lemma \ref{thmb_equiv}, we need an analogue of the limit functions $G_r$ from Theorem A for badly approximable irrationals. The following proposition shows that for large $k$, we can substitute the complicated expression of $P_{q_k}$ by a more controllable function $H_k$,
that is of similar shape to the representation of $G_{[k]}$ in \cite{quantum_invariants}.

\begin{prop}\label{limit_H}
For $h \in \mathbb{N}$, we define 
\[H_k(\alpha,\varepsilon) :=
2\pi \lvert \varepsilon + q_k\delta_k \rvert \prod_{n=1}^{\lfloor q_k/2\rfloor } h_{n,k}(\varepsilon),
\]
where
\[h_{n,k}(\varepsilon) = h_{n,k}(\alpha,\varepsilon) := \Bigg\lvert\Bigg(1 - q_k\delta_k\frac{\left\{n\cev{\alpha}_k\right\} - \frac{1}{2}}{n}\Bigg)^2 - \frac{\left(\varepsilon + \frac{q_k\delta_k}{2}\right)^2}{n^2}\Bigg\rvert.\]
For any badly approximable $\alpha$ and any compact interval $I$, we have 
\[P_{q_k}(\alpha,\varepsilon) = H_k(\alpha,\varepsilon)\left(1 + \mathcal{O}\left(q_k^{-2/3}\log^{2/3}q_k\right)\right) + \mathcal{O}(q_k^{-2}),\quad \varepsilon \in I,\]
with the implied constant only depending on $\alpha$ and $I$. In particular, we have

\[\lim_{k \to \infty} \lvert P_{q_k}(\alpha,x) - H_k(\alpha,x) \rvert = 0,\]
with the convergence being locally uniform on $\mathbb{R}$.
\end{prop}

\begin{proof}
We follow the strategy of \cite[Theorem 4]{quantum_invariants} and \cite[Lemma 3]{hauke_extreme},
which both treat the case where $\alpha$ is a quadratic irrational. Thus, the details in the proof that are identical to the arguments applied there are omitted.
We define $f(x) := \lvert 2\sin(\pi x)\rvert$ and 
\[ 
p_{n,k}(\varepsilon) = p_{n,k}(\alpha,\varepsilon) := \left\lvert\frac{f^2\Big(\frac{n}{q_k} - \big(\{n\cev{\alpha}_k\} - \frac{1}{2}\big)\delta_k\Big) - f^2\big(\frac{2\varepsilon + q_k\delta_k}{2q_k}\big)}{f^2\big(\frac{n}{q_k}\big)}\right\rvert.
\]
Using trigonometric identities and $\alpha = p_k/q_k + (-1)^k\delta_k/q_k$, we obtain

\begin{equation*}
    P_{q_k}(\beta,\varepsilon) =
    f(\delta_k + \tfrac{\varepsilon}{q_k})\,q_k
    \prod_{ 0 < n < q_k/2} \left\lvert\frac{f^2\Big(\frac{n}{q_k} - \big(\{n\cev{\alpha}_k\} - \frac{1}{2}\big)\delta_k\Big) - f^2\big(\frac{2\varepsilon + q_k\delta_k}{2q_k}\big)}{f^2\big(\frac{n}{q_k}\big)}\right\rvert
\end{equation*}
with an additional factor 
\begin{equation*}\frac{f\left(\frac{1}{2} - \frac{2\varepsilon + q_k\delta_k}{2q_k}\right)}{f(\frac{1}{2})} 
= 1 + \mathcal{O}(q_k^{-2})\end{equation*}
if $q_k$ is even.
Let $\psi(t) = t^{2/3}\log^{1/3} t$.
Since $\cev{\alpha}_k$ has bounded partial quotients, discrepancy estimates on
$\left\{ \{n\cev{\alpha}_k\}: 1 \ls n \ls N\right\}$ for $N < q_k$, together with 
Koksma's inequality (for details, see e.g. \cite{kuipers}) lead to 
\[
\prod_{\psi(q_k) < n < \frac{q_k}{2}} p_{n,k}(\varepsilon) = 1 + \mathcal{O}\Big(\frac{\log q_k}{\psi(q_k)}\Big), \quad 
\prod_{n = \psi(q_k)+1}^{\lfloor q_k/2 \rfloor} h_{n,k}(\varepsilon) = 1 + \mathcal{O}\Big(\frac{\log q_k}{\psi(q_k)}\Big).
\]
Furthermore, we can prove (see \cite{hauke_extreme}) that
\[
\prod\limits_{n=1}^{\psi(q_k)} p_{n,k}(\varepsilon) =
    \Big(1 + \mathcal{O}\Big(\frac{\psi(q_k)^2}{q_k^2}\Big)\Big){\prod\limits_{n=1}^{\psi(q_k)} h_{n,k}(\varepsilon)}
    + \mathcal{O}(q_k^{-2}),\]
 and clearly, we have
        \begin{equation*}
        f\Big(\delta_k + \frac{\varepsilon}{q_k}\Big)q_k =  2\pi\Big\lvert \varepsilon + q_k\delta_k\Big\rvert\Big(1 + \mathcal{O}\big(q_k^{-2}\big)\Big)
    + \mathcal{O}(q_k^{-2}).\end{equation*}
    Combining the previous estimates, we obtain the desired result.
\end{proof}

\begin{lem}\label{new_cond_badly}
Let $\alpha = [0;a_1,a_2,\ldots]$ be a badly approximable number, $a_K := \limsup_{n \to \infty} a_n \gs 2$ and $a_{\max} := \max_{n \in \mathbb{N}}a_n$.
Let $0 \ls x < y \ls 1, T \in \mathbb{N}$ and let 
\begin{equation*}
\begin{split}g(\ell,x,y) &:= \sum_{n =1}^{\ell} \frac{1}{2} - \left(\{nx\}\cdot\mathds{1}_{[\lfloor nx \rfloor = \lfloor ny \rfloor]}\right), \quad \\b(x) &:= \pi x\log x,\\
F(T,x,y) &:= \sum_{\ell = 1}^{T}\frac{g(\ell,x,y)}{\ell(\ell +1)}, \quad
\\E(T,a) &:= \frac{1+\log T}{T}\left(\frac{a}{8 \log a} + 6\right) + \frac{\frac{a}{8} + \frac{23}{4}}{T}.
\end{split}\end{equation*}
Assume there exists $\varepsilon'> 0, \in \mathbb{N}$ such that for infinitely many $k \in \mathbb{N}$ and $\cev{\alpha}_k \in [x_k,y_k]$, it holds that

\begin{equation}\label{cond_ineq}
b\left(\frac{a_{k+1}  + [0;\overline{a_K,1}] + x_k}{2\pi}\right) - F(T,x_k,y_k) - E(T,a_K) > \varepsilon'.
\end{equation}
Then we have
\[\liminf_{N \to \infty} P_N(\alpha) = 0, \quad \limsup_{N \to \infty} \frac{P_N(\alpha)}{N} = \infty.\]
\end{lem}
In order to prove this result, we need a technical estimate on the Birkhoff sum $\sum_{n = 1}^{\ell} f(n\alpha)$ where $f(x) := 1/2 -\{x\}$.
 \begin{prop}\label{pin_prop}
 Let $\alpha = [0;a_1,\ldots,a_k]$ be a rational number, $m \ls k$, $a := \max_{1 \ls i \ls m}a_i$ and let $\ell \in \mathbb{N}$ with $\ell < q_m$.
 Then we have
 
 \begin{equation*}
 \left\lvert \sum\limits_{n=1}^{\ell}\frac{1}{2} -\left\{n\alpha\right\}\right\rvert \ls  \left(\frac{a}{8 \log a} + 6\right)\log \ell + \frac{a}{8} + \frac{23}{4}.\end{equation*}
 \end{prop}
 
 \begin{proof}
 This is essentially \cite[Lemma 3]{grepstadII}, which is an immediate Corollary of \cite[Corollary 3]{pinner}. Note that in contrast to their setting, we deal with rational numbers, however, the proof remains the same as long as $\ell < q_m$.
 \end{proof}

\begin{proof}[Proof of Lemma \ref{new_cond_badly}]
We adopt the notations introduced in Proposition \ref{limit_H}.
Observe that if $\lvert \varepsilon\rvert$ is sufficiently small, we can remove the absolute values in the definition of $h_{n,k}(\alpha, \varepsilon)$. Since $\alpha$ is badly approximable, we have $q_k\delta_k \gg_{\alpha} 1$, so we can also assume $\lvert \varepsilon + q_k\delta_k \rvert > 0$. This leads to
\[h_{n,k}(\alpha,\varepsilon) \ls
1 + 2q_k\delta_k \frac{\frac{1}{2}- \{n \cev{\alpha}_k\}}{n} + \frac{2\varepsilon}{n^2}
\]
for any $n,k \in \mathbb{N}$, provided that $\varepsilon$ is sufficiently small. 
By $\log(1+x) \ls x$
and using again $q_k\delta_k \gg_{\alpha} 1$, we obtain

\begin{equation}\label{H_upperbound}
\begin{split}H_k(\alpha,\varepsilon) 
&\ls 2\pi (\varepsilon + q_k\delta_k)\cdot \exp\left(\sum_{n = 1}^{\lfloor q_k/2\rfloor}2q_k\delta_k \frac{\frac{1}{2}- \{n \cev{\alpha}_k\}}{n} + \frac{2\varepsilon}{n^2}\right)
\\&= \left(1 + \mathcal{O}_{\alpha}(\delta)\right)\cdot 2\pi q_k\delta_k\exp\left(2q_k\delta_k\sum_{n = 1}^{\lfloor q_k/2\rfloor} \frac{\frac{1}{2}- \{n \cev{\alpha}_k\}}{n}\right).
\end{split}
\end{equation}
Now let $T \in \mathbb{N}$ with $T < q_k/2$.
  Writing 
     $S_{\ell}(\alpha) :=  \sum\limits_{n=1}^{\ell}\frac{1}{2} -\left\{n\alpha\right\}$, we use summation by parts to obtain
     
     \begin{equation*}
     \begin{split}\sum_{n=1}^{\lfloor q_k/2\rfloor}\frac{\frac{1}{2} -\left\{n\cev{\alpha}_k\right\}}{n}
&\ls \sum_{\ell=1}^{T}\frac{S_{\ell}(\cev{\alpha}_k)}{\ell(\ell +1)} +\frac{S_{\lfloor q_k/2\rfloor}}{\lfloor q_k/2\rfloor}
 + \sum_{\ell=T+1}^{\lfloor \sqrt{q_k}\rfloor}\frac{\lvert S_{\ell}(\cev{\alpha}_k)\rvert}{\ell^2}
 + \sum_{\ell=\lfloor \sqrt{q_k}\rfloor+1}^{\lfloor q_k/2\rfloor}\frac{\lvert S_{\ell}(\cev{\alpha}_k)\rvert}{\ell^2}.
 \end{split}
 \end{equation*}
 We apply Proposition \ref{pin_prop} to the last three terms on the right-hand side of the previous inequality. 
 Since $\limsup_{i \to \infty} a_i = a_K$, there exists some $K_0 \in \mathbb{N}$ such that 
 $\max_{k \gs K_0} a_i = a_K$.
 If $k$ is sufficiently large, we have
 that $\sqrt{q_k} < q_{k-K_0}(\alpha_k)$
 (this follows immediately from the fact that $q_k(\alpha_k) = q_k(\alpha)$). Hence, we obtain
 \[\sum_{\ell=T+1}^{\lfloor \sqrt{q_k}\rfloor}\frac{\lvert S_{\ell}(\cev{\alpha}_k)\rvert}{\ell^2}
 \ls E(T,a_K),
 \]
 and another application of Proposition \ref{pin_prop} with $a = a_{\max}$ yields
 
 \[\frac{S_{\lfloor q_k/2\rfloor}}{\lfloor q_k/2\rfloor}  + \sum_{\ell=\sqrt{q_k}}^{\lfloor q_k/2\rfloor}\frac{\lvert S_{\ell}(\cev{\alpha}_k)\rvert}{\ell^2}
 = \mathcal{O}_{\alpha}\left(\frac{\log q_k}{\sqrt{q_k}}\right).
 \]
 By a short case distinction argument, we see that $S_{\ell}(\cev{\alpha}_k) \ls g(\ell,x,y)$ holds for all $\ell \in \mathbb{N}$, hence combined with the previous estimates we obtain
\begin{equation*}
\sum_{n=1}^{\lfloor q_k/2 \rfloor}\frac{\frac{1}{2} -\left\{n\cev{\alpha}_k\right\}}{n} \ls F(T,x,y) +
E(T,a_K) + \mathcal{O}_{\alpha}\left(\frac{\log q_k}{\sqrt{q_k}}\right), \quad k \to \infty.\end{equation*}
Applying a logarithm to \eqref{H_upperbound} leads to
\begin{equation}\begin{split}\label{Hlog_upperbound}
    \frac{\log H_k(\alpha,\varepsilon)}{2q_k\delta_k} &\ls 
    \mathcal{O}_{\alpha}(\delta) + F(T,x,y) +
E(T,a_K) + \mathcal{O}_{\alpha}\left(\frac{\log q_k}{\sqrt{q_k}}\right)
- b\left(\frac{1}{2\pi q_k\delta_k}\right).
\end{split}
\end{equation}
Observe that for $K \gs K_0$, we have
\[\frac{1}{q_k\delta_k} = a_{k+1} + \cev{\alpha}_k + \vec{\alpha}_{k+2} 
\gs 
a_{k+1}  + [0;\overline{a_K,1}] + x_k,\]
where we used that $\vec{\alpha}_{k+2} \in \mathcal{Q}(a_K)$ and $\cev{\alpha}_k \in [x_k,y_k]$.
Since $b$ is monotonically increasing on $[0,\infty)$, this implies
\begin{equation*}
   b\left(\frac{a_{k+1}  + [0;\overline{a_K,1}] + x_k}{2\pi}\right) \ls  b\left(\frac{1}{2\pi q_k\delta_k}\right).
\end{equation*}
Using the assumption \eqref{cond_ineq} and choosing $k, \delta$ in \eqref{Hlog_upperbound}
sufficiently large respectively sufficiently small, we deduce that for infinitely many $k$ and any $\varepsilon \in (-\delta,\delta)$, we have
$\log H_k(\alpha,\varepsilon) < -\varepsilon'/2$. Using Proposition \ref{limit_H}, this shows that for $k$ large enough, we fulfill \eqref{shifted<1} and thus, the statement follows from Lemma \ref{thmb_equiv}.
\end{proof}

\begin{figure}[h]
  \centering
  \includegraphics[width=.7\linewidth]{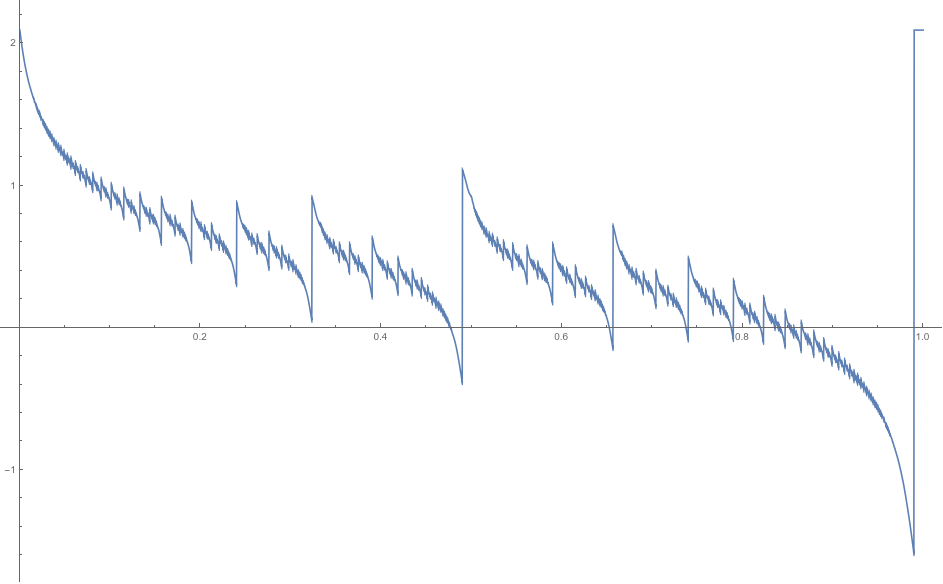}
  \caption{Values of $\sum_{\ell = 1}^{100}\frac{g(\ell,x_j,y_j)}{\ell(\ell +1)}$ for $y_j = x_j + 1/100$.
  We see that the global maximum appears at $x = 0$, with local peaks next at rationals with small denominator.
  }
  \label{fig:sfig1}
\end{figure}

Now we are in position to prove Theorem \ref{K_l<=7} and start with a short outline of the proof.
Looking closer at the condition \eqref{cond_ineq} in Lemma \ref{new_cond_badly}, and assuming $T$ to be very large and $x_k,y_k$ very close together, this is morally equivalent to

\begin{equation}\label{morally}F(T,x_k,y_k)
< b\left(\frac{a_{k+1}  + [0;\overline{a_K,1}] + x_k}{2\pi}\right)
< \frac{a_{k+1}}{2} \log\left(\frac{a_{k+1}}{2\pi}\right)\end{equation}
when $\cev{\alpha}_k \in [x_k,y_k]$.
For large $a_K$, we can use the trivial estimate $g(\ell,x,y) \ls \ell/2$ to show that
\eqref{morally} holds for any $x \in [0,1]$, so we can concentrate on reasonably small $a_K$.
We see in Figure \ref{fig:sfig1} that $F(T,x,y)$ is not too big when bounded away from $0$ and (other) rationals with a small denominators. If $\alpha$ is a badly approximable number, then
$\cev{\alpha}_k \in \mathcal{B}(a,M)$ with $M$ large is bounded away from these rationals (see Proposition \ref{not_close_to_rational} below).
Covering $\mathcal{B}(a,M)$ with a finite number of small intervals $[x,y]$ that avoid being close to rationals with small denominators, we are left to prove that \eqref{cond_ineq} respectively \eqref{morally} holds on each such interval, which reduces the problem to checking finitely many cases. As we are bounded away from the values $x,y$ where $F(T,x,y)$ is large, we can morally think of $F(T,x,y) \ls 0$, so the question is essentially reduced to whether $a_{k+1} > 2\pi \approx 6.28$. 
A refinement of this argument is in fact enough to prove Theorem \ref{K_l<=7} by considering those $k$ where $a_{k+1} = a_K \gs 8$, and does also work for most badly approximable numbers with $a_K = 7$. However, it might happen that $a_{k+1} = 7$ and $\cev{\alpha}_{k} \approx x_k$ is so small that $F(T,x,y)$ is significantly larger than $0$, and our estimates are too coarse to show \eqref{cond_ineq} directly. Fortunately, this only happens when $a_{k} = 6$ and $\cev{\alpha}_{k-1}$ is so large that $F(T,x,y)$ for $x \approx \cev{\alpha}_{k-1}$ is significantly smaller than $0$, enabling us to show \eqref{cond_ineq} for $k-1$ instead.

\begin{prop}\label{not_close_to_rational}
Let $a \gs 2, M \gs 100, 1 \ls m \ls a+1$.
Then $\mathcal{B}(a,M)$ is disjoint from the following sets:
\begin{align}\label{disj1}&\left[0,\frac{1}{a+1}\right], \quad \left[\frac{a+1}{a+2},1\right], \quad \\&\label{disj2}
\left[\frac{1}{m} - \frac{1}{m^2(a+2)}, \frac{1}{m} + \frac{1}{m^2(a+2)}\right],
\\&
\left[\frac{2}{2m+1}- \frac{1}{(2m+1)^2(a+3)},\frac{2}{2m+1}+ \frac{1}{(2m+1)^2(a+3)}\right],\label{disj3}
\\& 
\left[\frac{3}{3m+1}- \frac{1}{(3m+1)^2(a+3)},\frac{3}{3m+1}+ \frac{1}{(3m+1)^2(a+3)}\right],\label{disj4}
\\& 
\left[\frac{3}{3m+2}- \frac{1}{(3m+2)^3(a+3)},\frac{3}{3m+2}+ \frac{1}{(3m+2)^2(a+3)}\right]\label{disj5}
.\end{align}
\end{prop}

\begin{proof}
Clearly, $\alpha \in \mathcal{B}(a,M)$ fulfills $\alpha > [0;a,1] = \frac{1}{a+1}$
and $\alpha < [0;1,a+1] = \frac{a+1}{a+2}$, hence $\mathcal{B}(a,M)$ is disjoint from the first two sets in \eqref{disj1}.
Now let $\alpha \in \mathcal{B}(a,M) \cap [\frac{1}{m} - \frac{1}{m^2(a+2)},\frac{1}{m}]$, 
then $\alpha = [0;m,\alpha']$ for some $\alpha' \in \mathcal{B}(a,M-1)$. By \eqref{disj1} (which still holds for $M-1$), we have $\alpha' > \frac{1}{a+1}$, hence $\alpha < \frac{1}{m} - \frac{1}{m^2(a+1 + 1/m)} \ls \frac{1}{m} - \frac{1}{m^2(a+2)}$, a contradiction. Similarly, if
 $\alpha \in \mathcal{B}(a,M) \cap [\frac{1}{m},\frac{1}{m} + \frac{1}{m^2(a+2)}]$, then
 $\alpha = [0;m+1,\alpha']$ and again by \eqref{disj1} $\alpha' < \frac{a+1}{a+2}$, we have 
 $\alpha > \frac{1}{m} + \frac{1}{m^2(a+2)}$.
 Now assume that $\alpha \in \mathcal{B}(a,M) \cap \left[\frac{2}{2m+1}- \frac{1}{(a+3)(2m+1)^2}, \frac{2}{2m+1} + \frac{1}{(a+3)(2m+1)^2}\right]$. Note that 
 \eqref{disj3} is contained in $\left(\frac{1}{m},\frac{1}{m+1}\right)$, hence
 $\alpha = \frac{1}{m + \alpha'}$ with $\alpha' \in \mathcal{B}(a,M)$. 
   If $\alpha' \gs \frac{1}{2}$, then using
 \eqref{disj2} with $m = 2$ implies
 $\alpha' > \frac{1}{2}+ \frac{1}{4(a+2)}$ and thus, we have
 \begin{equation*}
     \begin{split}
 \alpha < \frac{1}{m + \frac{1}{2} + \frac{1}{4(a+2)}}
 &= \frac{2}{2m+1} - \frac{2}{(2(a+2)(2m+1)+1)(2m+1)}
 \\&< \frac{2}{2m+1} - \frac{1}{(a+3)(2m+1)^2}.
 \end{split}
 \end{equation*}
 If $\alpha' < \frac{1}{2}$, then again by \eqref{disj2}, $\alpha' < \frac{1}{2} - \frac{1}{4(M+2)}$, and so
  \begin{equation*}
     \begin{split}
 \alpha > \frac{1}{m + \frac{1}{2} - \frac{1}{4(a+2)}}
 &= \frac{2}{2m+1} + \frac{2}{(2(a+2)(2m+1)-1)(2m+1)}
 \\&> \frac{2}{2m+1} + \frac{1}{(a+3)(2m+1)^2}.
 \end{split}
 \end{equation*}
 The proof of \eqref{disj4} and \eqref{disj5} works precisely in the same fashion as \eqref{disj3}, by applying \eqref{disj2} with $m =3$ respectively $m =2$.
 \end{proof}

\begin{proof}[Proof of Theorem \ref{K_l<=7}]
Let $\alpha = [a_0;a_1,a_2,\ldots]$ be a badly approximable irrational 
(the case where $\alpha$ has unbounded partial quotients was proven by Lubinsky in \cite{lubinsky})
with $a_K = \limsup_{k \to \infty} a_k \gs 7$ and let $K_0 \in \mathbb{N}$ such that $\alpha \in \mathcal{B}(a_K,K_0)$.
If $\limsup_{k \to \infty}a_k = \liminf_{k \to \infty}a_k$, then there exists some $M \in \mathbb{N}$ such that
$T^{M}(\alpha) = [0;\overline{a_K}]$ (with $T$ denoting the Gauss map).
It was shown in \cite{tech_zaf} that 
$\liminf_{N \to \infty} P_N([0;\overline{a_K}]) = 0$ and arguing as in \cite{grepstadII}, this behaviour is invariant under the application of the Gauss map.
Thus, we will assume from now on that $\liminf_{k \to \infty}a_k < a_K-1$. This implies that there exists an increasing sequence
$(k_j)_{j \in \mathbb{N}}$ that fulfills the following for every $j \in \mathbb{N}$:

\begin{itemize}
    \item $k_j \gs K_0 +1000$.
    \item $a_{k_j+1} = a_K$.
    \item $a_{k_j} \ls a_K-1$.
\end{itemize}
We distinguish cases for the value of $a_K$.\vspace{3mm}

\begin{itemize}
\item $a_K \gs 18$:
In view of Lemma \ref{new_cond_badly}, it suffices to show \eqref{cond_ineq} for all $k_j$ and some $T \in \mathbb{N}$. Using the trivial bound $g(\ell,x,y) \ls \frac{1}{2}$ and monotonicity of $b$, we have
\[F(T,x,y) < \frac{\log(T+1)}{2}, \quad b\left(\frac{a_K  + [0;\overline{a_K,1}] + x}{2\pi}\right)
> \frac{a_{K}}{2} \log\left(\frac{a_{K}}{2\pi}\right),\] so we are left to prove that

\begin{equation}\label{case_large_ak}
\frac{\log(T+1)}{a_K} + \frac{2E(T,a_K)}{a_K} < \log\left(\frac{a_K}{2\pi}\right) + \varepsilon',
\end{equation}
which is fulfilled for $T = 50$ and $a_K = 18$.
Since the left-hand side of \eqref{case_large_ak} is monotonically decreasing in $a_K$, whereas the right-hand side is increasing, the result follows immediately for all $a_K \gs 18$. 

\item $8 \ls a_K \ls 18$:
As in the case $a_K \gs 18$ we fix an arbitrary $k_j$, but instead of working with trivial bounds on $g(\ell,x,y)$, we compute the sum explicitly in small grids for $[x,y]$. 
As $E(T,a_K)$ and $b\left(\frac{a_K  + [0;\overline{18,1}]+x}{2\pi}\right)$
are monotonically increasing respectively decreasing in $a_K$, \eqref{cond_ineq} can be deduced from
\[F(x,y,T) + E(T,18) < b\left(\frac{9  + [0;\overline{18,1}]+x}{2\pi}\right)\]
whenever $9 \ls a_K \ls 18$.
We can prove with computational assistance that for $T = 200, R = 2000$, we have
\begin{equation*}
	F\left(\frac{i}{R},\frac{i+1}{R},T\right) + E(T,18) < b\left(\frac{9  + [0;\overline{18,1}]+x}{2\pi}\right)\end{equation*}
for all integers $i$ that satisfy $\left\lfloor \frac{R}{19} \right\rfloor \ls i \ls \left\lceil \frac{20 R}{21}\right\rceil$.
By Proposition \ref{not_close_to_rational}, $\cev{\alpha}_{k_j-1} \in \mathcal{B}(18,10) \subseteq [1/19,19/20]$.
Hence, \eqref{cond_ineq} holds for any badly approximable $\alpha$ with
$9 \ls a_K \ls 18$ and thus, Lemma \ref{new_cond_badly}
gives the result.
Analogously, we can prove for the same values of $T,R$ that we have

\begin{equation*}F\left(\frac{i}{R},\frac{i+1}{R},T\right) + E(T,8) < b\left(\frac{8  + [0;\overline{8,1}]+x}{2\pi}\right)\end{equation*}
for all integers $i$ that satisfy $\left\lfloor \frac{R}{9} \right\rfloor \ls i \ls \left\lceil \frac{9R}{10}\right\rceil$ and can conclude the case $a_K = 8$ exactly in the same way as for $9 \ls a_K \ls 18$.
\item $a_K = 7$:
This case is by far the most intricate one.
We cut out small windows close to rationals with a small denominator, justified by Proposition \ref{not_close_to_rational}, in the following way: we lay a grid of size $R$ over the unit interval and remove those intervals $\left[\frac{i}{R},\frac{i+1}{R}\right]$ that are completely contained inside one of the sets \eqref{disj1}\,--\,\eqref{disj5} with $m = 1,\ldots,6$.
In that way we can prove using $T = 4.000, R = 360.000$ that \eqref{cond_ineq} holds whenever
\begin{equation*}
a_{k_j+1} = 7,\quad \cev{\alpha}_{k_j} \gs [0;6,1,6,1,7,1,7]
\approx 0.14549
 .\end{equation*}
\begin{figure}
\centering
\begin{subfigure}[b]{\textwidth}
  \centering
  \includegraphics[width=\linewidth]{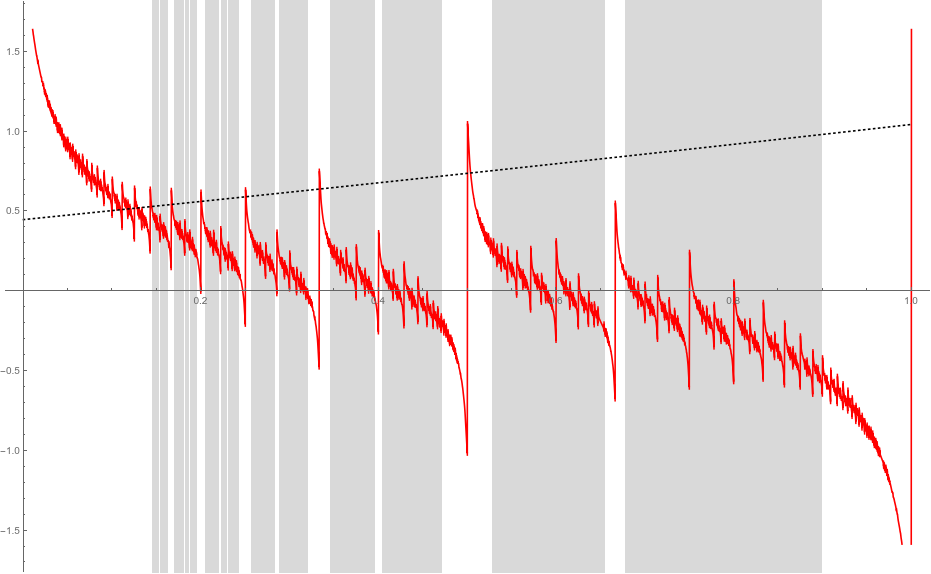}
  \caption*{The case $a_K = 7, \cev{\alpha}_{k_j} \gs [0;6,1,6,1,7,1,7]$.
  By making use of Proposition \ref{not_close_to_rational}, we are allowed to remove the peaks coming from rationals of the form $1/m, 2/(2m+1), 3/(3m+1), 3/(3m+2), m = 1,\ldots, 5,6$. We see that
  everywhere in the grey-shaded region, except when $\cev{\alpha}_{k_j}$ is very close to
  $[0;6,1,7]$, the inequality \eqref{cond_ineq} holds. 
  }
   \label{fig:caseak71} 
\end{subfigure}
\begin{subfigure}[b]{\textwidth}
  \centering
  \includegraphics[width=\linewidth]{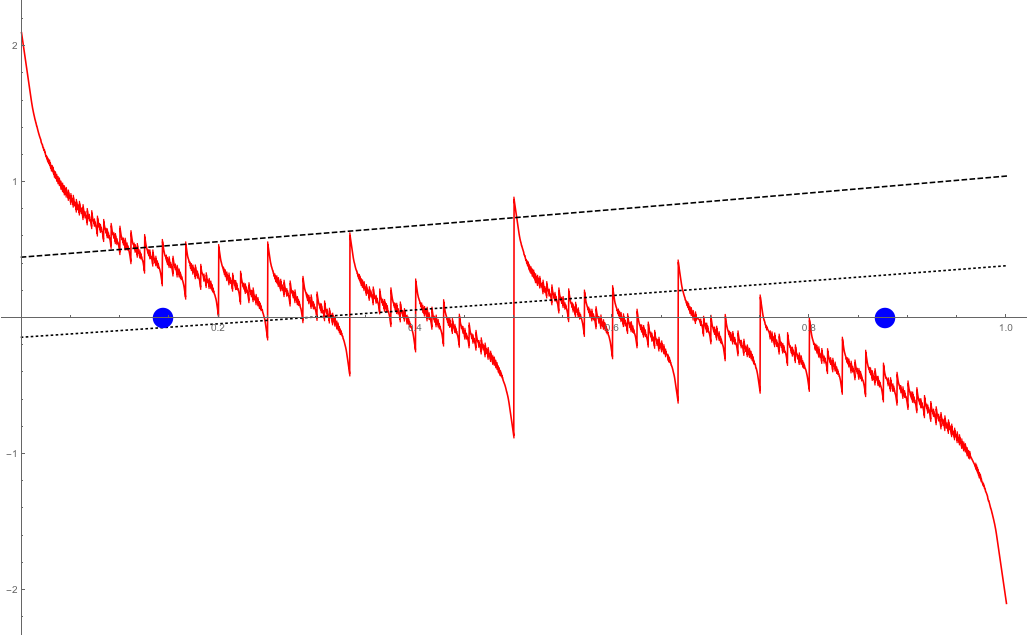}
  \caption*{The case $a_K = 7, \cev{\alpha}_{k_j} < [0;6,1,6,1,7,1,7]$,
  illustrated at the irrational $\alpha = [0;\overline{1,6,7}]$. The dashed line is the bound $b\left(\frac{7  + [0;\overline{7,1}]+x}{2\pi}\right)$, whereas
  the dotted line below is $b\left(\frac{6  + [0;\overline{7,1}]+x}{2\pi}\right)$.
  The big dot close to zero represents $\cev{\alpha}_{3m +2}$ ($m$ a large integer), where we cannot prove that
  \eqref{cond_ineq} holds. However, this forces $\cev{\alpha}_{3m+1}$ (represented by the big dot close to $1$) to be so large that \eqref{cond_ineq} holds also for $a_{3m +2} = 6$.}
   \label{fig:caseak72}
\end{subfigure}
\end{figure}
This proves Theorem \ref{K_l<=7} for badly approximable $\alpha$ with $a_K = 7$
such that\\
$\cev{\alpha}_{k_j} \gs [0;6,1,6,1,7,1,7]$ holds infinitely often.
If this is not the case, using $a_{k_j-1} \ls 6$, this implies that
$a_{k_j} = 6, a_{k_j-1} = 1, a_{k_j-2} = 7$ holds infinitely often.
Therefore, we can deduce that
$\cev{\alpha}_{k_j-1} \gs [0;1,7,7,1,7,1] > 0.876$.
For $T = 2.000,R = 50.000$, we have that
\begin{equation*}
F\left(\frac{i}{R},\frac{i+1}{R},T\right) + E(T,7) < b\left(\frac{6  + [0;\overline{7,1}]+x}{2\pi}\right)\end{equation*}
holds for all integers $i$ that satisfy $\left\lfloor 0.87\cdot R \right\rfloor \ls i \ls \left\lceil \frac{9R}{10}\right\rceil$.
Equation \eqref{cond_ineq} follows now by the same argumentation as above.
\end{itemize}
\end{proof}

\section{Proof of Theorem \ref{l=2}}\label{secl2}
From now on, we will only deal with quadratic irrationals of the form\\$\alpha = [0;a_1,\ldots,a_p,\overline{a_{p+1},\ldots,a_{\ell}}]$ with $a_K = \max_{p+1 \ls i \ls p+{\ell}} a_i$. By invariance under application of
 the Gauss map (which will follow from Lemma \ref{aist_borda_limit_fct} and the method of proof naturally), we can assume that $\alpha \in \mathcal{Q}(a_K)$. As the statement of Theorem \ref{l=2} follows for $a_K \gs 7$ from Theorem $\ref{K_l<=7}$,
we can fix $a_K = 6$, and thus only have to check the five quadratic irrationals $\alpha = [0;\overline{6,a}]$ with $a \in \{1,\ldots,5\}$.
For each fixed quadratic irrational, the limit functions $G_r(\alpha,\cdot), r \in \{0,1\}$ from Theorem A can be computed explicitly, up to an arbitrarily small error
(see Lemma \ref{approx_limit} below).
If $\alpha = [0;\overline{6,a}]$
 with $a \in \{1,\ldots,4\}$, we will show that
$G_{0}(\alpha,0) < 1$ and thus, the statement follows immediately from Theorem B.
If $\alpha = [0;\overline{6,5}]$, Theorem B fails since
$\min \{ G_{0}(\alpha,0),G_{1}(\alpha,0)\} > 1$. However, we adapt the proof strategy of Aistleitner, Technau  \cite{tech_zaf} that was used in order to show that
$\alpha = [0;\overline{6}]$ fulfills \eqref{asym=0}: instead of choosing 
an Ostrowski expansion with only very few non-zero coefficients,
we consider an Ostrowski expansion with all coefficients being $1$. For $q_{k} \ls N < q_k$, this leads to almost all negative perturbations fulfilling $\varepsilon_i(N) \approx -0.025, 1 \ls i \ls k$ and so, $P_{q_{i+1}}(\alpha,\varepsilon_{i+1}(N)) \cdot P_{q_{i}}(\alpha,\varepsilon_i(N)) \lessapprox 0.999 < 1$
holds for most such $i$. Using the decomposition \eqref{shifted_prod_1}, we will deduce the result.To formalize this, we start with the ingredients needed to explicitly compute the value of
$G_r(\alpha,\varepsilon)$ for fixed $\alpha,\varepsilon$.
We use a representation of the limit functions from Theorem A that is due to Aistleitner and Borda \cite{quantum_invariants}, which is much easier to control than the one used in \cite{grepstadII} (although by
the uniqueness of limits, the functions obviously coincide).

\begin{lem}[see \cite{quantum_invariants}, Theorem 4]\label{aist_borda_limit_fct}
Let $\alpha = [0;a_1,\ldots,a_p,\overline{a_{p+1},\ldots,a_{\ell}}]$ and let $\alpha_{\tau_r}$ and $C(r)$ as in \eqref{alphatau} and \eqref{quadr_limit1}. Further, define
\begin{equation*}
    G_{r}(\alpha,\varepsilon) = 2\pi\lvert \varepsilon + C(r)\rvert\prod_{n=1}^{\infty}
     \left\lvert g_{n,r}(\alpha_{\tau_r},\varepsilon) \right\rvert
\end{equation*}
where
\begin{equation*}
g_{n}(\alpha_{\tau_r},\varepsilon) := \left(1 - C(r)\frac{\left\{n\alpha_{\tau_r}\right\} - \frac{1}{2}}{n}\right)^2 - \frac{\left(\varepsilon + \frac{C(r)}{2}\right)^2}{n^2}.\end{equation*}
Then we have
\begin{equation}\label{convergence_rate}
    \lim_{m \to \infty}P_{q_{m\ell + r+p}}(\alpha,\varepsilon) = G_{r}(\alpha,\varepsilon).
\end{equation}
\end{lem}

\begin{rmk}
In case of quadratic irrationals, one can show that the functions $H_k$ defined in Proposition \ref{limit_H} converge along the subsequences $(m \ell + r)_{m \in \mathbb{N}}, 0 \ls r \ls \ell -1$ locally uniformly to $G_r$. The convergence rate in \eqref{convergence_rate} is also the same as the one obtained in Proposition \ref{limit_H}.
\end{rmk}

\begin{lem}\label{approx_limit}\item
\begin{enumerate}[(i)]
\item Let $[a,b] \subseteq I$ where $I$ is a zero-free interval of $G_r(\alpha,\cdot)$. 
Then we have
\begin{equation*}
\min\left\{G(a),G(b)\right\} \ls G(x) \quad \forall x \in [a,b].\end{equation*}
If $[a,b]$ does not contain the maximizer (which is unique in $I$) of $G_r(\alpha,\cdot)$, then we also have
\begin{equation*}
G(x) \ls \max\left\{G(a),G(b)\right\}\quad  \forall x \in [a,b].\end{equation*}
\item Let $T \in \mathbb{N}, 0 \ls r \ls K-1$ and $\varepsilon \in I \subseteq (-1,1) $ 
such that $G_r(\varepsilon) > 0$ for all $\varepsilon \in I$.
Let
$G_{r,T}(\alpha,\varepsilon) := 2\pi\lvert \varepsilon + C(r)\rvert\prod_{n=1}^{T}
     \left\lvert g_{n,r}(\cev{\alpha}_r,\varepsilon)\right\rvert$  
     Then we have for $T$ sufficiently large that
     
    \begin{align}\label{G_T_ineq}
     \left(1 - 2C(r)E(T,a_K) - \frac{5}{T}\right)G_{r,T}(\varepsilon) &\ls G_{r}(\varepsilon)
     \ls G_{r,T}(\varepsilon)\left(1 - 3C(r)E(T,a_K)\right)^{-1}
     .\end{align}
\end{enumerate}
\end{lem}

\begin{proof}
\begin{enumerate}[(i)]
          \item 
        This follows directly from the fact that $G_{r}(\alpha,\varepsilon)$ is log-concave on zero-free intervals, which can be proven in the precisely same way as it was done in 
        \cite[Proposition 1]{tech_zaf} for $\alpha$ being the golden ratio.
    \item 
For the second inequality in \eqref{G_T_ineq}, we follow the strategy of \cite[Lemma 3]{hauke_extreme}. Clearly, it suffices to show that

\begin{equation}\label{suffices_limit}\prod_{n=T+1}^{\infty}
     \left\lvert g_{n,r}(\cev{\alpha}_r,\varepsilon)\right\lvert \ls 1 + 3C(r)E(T,a_K).\end{equation}
     Observe that 
     $g_{n,r}(\cev{\alpha}_r,\varepsilon) \gs 1 - 2/n - 4/n^2,$
     hence we can remove the absolute values around $g_{n,r}$ in the definition of $G_{r}$ for $n \gs 4$.
     Furthermore, this implies
   
     \[\prod_{n=T+1}^{\infty}
     g_{n,r}(\cev{\alpha}_r,\varepsilon) \ls \prod_{n=T+1}^{\infty}\left(1 - C(r)\frac{\left\{n\cev{\alpha}_r\right\} - \frac{1}{2}}{n}\right)^2.\]
     By arguments as in Lemma \ref{new_cond_badly}, we obtain
     \begin{equation}\label{pinner_applied}\sum_{n=T+1}^{\infty}
     C(r)\frac{\left\{n\cev{\alpha}_r\right\} - \frac{1}{2}}{n}
     \ls C(r)E(T,a_K),
     \end{equation}
     so we can deduce  
     \begin{equation*}
         \begin{split}
     \prod_{n=T+1}^{\infty}
     g_{n,r}(\cev{\alpha}_r,\varepsilon)
     &\ls \prod_{n=T+1}^{\infty}\left(1 - 2C(r)\frac{\left\{n\cev{\alpha}_r\right\} - \frac{1}{2}}{n}\right)
     \ls
     \exp\left(\sum_{n=T+1}^{\infty}- 2C(r)\frac{\left\{n\cev{\alpha}_r\right\} - \frac{1}{2}}{n}\right)
     \\&\ls
     \exp\left(2C(r)E(T,a_K)\right).
     \end{split}
        \end{equation*}
          For $T$ sufficiently large, $2C(r)\cdot E(T,a_K)$ is bounded from above by $1/2$, and for $x < \frac{1}{2},$ we have the inequality $\exp(x) \ls 1 + \frac{3}{2}x$, which leads to \eqref{suffices_limit}.\par
          
          Concerning the first inequality in \eqref{G_T_ineq}, we use that
          
          \[\prod_{n = N}^M \left(1 + a_n\right) \gs 1 - \left\lvert \sum_{n = N}^M a_n \right\rvert 
- \frac{1}{N-1}\]
          holds for any $N, M \in \mathbb{N}$, $\lvert a_n \rvert \ls \min\{\frac{1}{2}, 1/n\}$, a fact that can be proven by simple estimates on the Taylor series of $\log$ and $\exp$ 
(see  \cite[Lemma 9]{hauke_extreme}). This implies for $T$ sufficiently large that
          
          \[\begin{split}\prod_{n = T+1}^{\infty} g_{n,r}(\cev{\alpha}_r,\varepsilon) \gs
          1 &- 2C(r)\left\lvert\sum_{n = T+1}^{\infty}  \left(\frac{\{n \cev{\alpha}_r\} - \frac{1}{2}}{n}\right)\right\rvert \\&- \sum_{n = T+1}^{\infty}\frac{C(r)^2\left(\{n \cev{\alpha}_r\} - \frac{1}{2}\right)^2 + \left(\varepsilon + C(r)/2\right)^2}{n^2}
          - \frac{1}{T}.\end{split}
          \]
          Applying \eqref{pinner_applied} and $\lvert C(r)\rvert , \lvert \varepsilon \rvert < 1$, we obtain
          
          \[\prod_{n = T+1}^{\infty} g_{n,r}(\cev{\alpha}_r,\varepsilon)
          \gs 1 - 2C(r)E(T,a_K) - \frac{5}{T}.
          \]

       \end{enumerate}
\end{proof}

\begin{rmk}\label{explicit_com}
Using Lemma \ref{approx_limit}, we can now compute lower respectively upper bounds for $G_r(\alpha,\varepsilon)$ that approximate the exact value arbitrarily well, by taking $T$ large and small intervals $[a,b] \in \varepsilon$ 
such that $a,b$ can be represented in a computer program symbolically (e.g. rational numbers).
In the rest of the paper, whenever we compute a value $G_r(\alpha,\varepsilon)$, we will always implicitly apply this procedure, without stating the exact intervals and the size of $T$.
Note that with this approach, we can also prove the following: if $\alpha = [0;\overline{1,a_2}]$,  then $G_1(\alpha,0) < 1$ if and only if $a_2 \gs 4$, and if $\alpha = [0;\overline{2,a_2}]$, then  $G_1(\alpha,0) < 1$ if and only if $a_2 \gs 5$, statements which have been conjectured in \cite{grepstadII}.
\end{rmk}

In order to prove $\liminf_{N \to \infty} P_N\left([0;\overline{6,5}]\right) = 0$, we need to study the behaviour of $P_{q_k}(\alpha,\varepsilon)$ when $\varepsilon$ is bounded away from $0$. To do so, we need to define some more notations.
We call a tuple $(\beta_1,\ldots,\beta_j)$ \textit{admissible with respect to $r$ and $\alpha$} 
(where $0 \ls r \ls \ell-1$) if there exists some $N \in \mathbb{N}, i = [r]$ 
such that if
$N = \sum_{j =1}^n b_j(N)q_j(\alpha)$ is the Ostrowski expansion of $N$, we have
$(b_{i-1},b_{i},\ldots,b_{i+j-1}) = (\beta_1,\beta_2,\ldots,\beta_j)$.
Similarly, we define a sequence
$(\beta_j)_{j \in \mathbb{N}}$ to be admissible if for any $J \in \mathbb{N}$, $(\beta_j)_{1 \ls j \ls J}$ is admissible with respect to $r = 0$ and $\alpha$. The purpose of this definitions is that the admissibility of a sequence respectively tuple
encodes the rules on the digits of the Ostrowski expansion. For example, we call a tuple $(a,b,c) \in \mathbb{N}^3$ admissible with respect to $\alpha = [0;\overline{6,5}]$ and $1$ if $a \leq 6, b \leq 5, c \leq 6$, if $b = 5$ then $a = 0$ and if $c = 6$, then $b = 0$.
For convenience, we will drop the dependence of the admissibility of $(\beta_1,\ldots,\beta_j)$ on $r$ and $\alpha$, whenever $\alpha,r$ are implicitly defined by the context. 
For $(\beta_j)_{j \in \mathbb{N}}$ being an admissible sequence, we define for $0 \ls r \ls \ell-1, [i] = r$ and $1 \ls k \ls \beta_i$,
\begin{equation}\label{def_eps'_seq}
\varepsilon_{r,k}'((\beta_{i+j})_{j \in \mathbb{N}})
:= kC(r) 
+ \sum_{j =1}^{\infty}(-1)^j \beta_{r+j}C([r+j])\cdot \prod_{n = 1}^{j}\alpha_{\sigma_{[r+n]}},\end{equation}
whenever the sum converges.
The motivation behind this definition is the following: 
suppose that $N_n$ is the sequence of integers along which we hope to prove that 
$\lim_{n \to \infty} P_{N_n}(\alpha) = 0$, and assume that
$N_n= \sum_{m =1}^n \beta_mq_m$ with $(\beta_j)_{j \in \mathbb{N}}$ being an admissible sequence.
If $(\beta_j)_{j \in \mathbb{N}}$ has a specific periodic structure (say for convenience, with same period $\ell$ than the continued fraction expansion of $\alpha$), then we can explicitly compute the actual value of $\varepsilon_{r,k}'((\beta_{i+j})_{j \in \mathbb{N}})'$,
which no longer depends on $i$, but only on $[i]=r$. The following proposition tells us that we can
replace the perturbation value of $\varepsilon_{i,k}(N_n)$ by $\varepsilon_{[i],k}'$
and thus, $P_{q_i}(\varepsilon_{i,k}(N_n)) \approx P_{q_i}(\varepsilon_{[i],k}')
\approx G_{[i]}(\varepsilon_{[i],k}')$, with the value of the latter term being explicitly computable, so we only need to compute $\ell$ function evaluations, regardless of $n$.

\begin{prop}\label{eps_convergence}

Let $(\beta_j)_{j \in \mathbb{N}}$ be an admissible sequence and let $(N_n)_{n \in \mathbb{N}}$ with
$N_n := \sum_{m =1}^n \beta_mq_m$ the associated sequence. For every $\delta > 0$, there exist $I_0,J_0$ such that for all $n \gs J_0$,
\begin{equation}\label{conv_of_eps}
    \sup_{I_0 < i < n - J_0}  \left\lvert\varepsilon_{i,k}(N_n) - \varepsilon'_{[i],k}\left((\beta_{i+j})_{j \in \mathbb{N}}\right) \right\rvert < \delta.
\end{equation}
\end{prop}

\begin{proof}
If $i \ls n - J_0$, we have
\begin{equation}\label{I0J0}\begin{split}\lvert \varepsilon_{i,k}(N_n)
- \varepsilon'_{r,k}((\beta_{i+j})_{j \in \mathbb{N}})\rvert
&\ls 
\sum_{j = 1}^{J_0}
\beta_{i+j}\left\lvert (-1)^iq_i\delta_{i+j} - (-1)^jC([r+j])\prod_{n=1}^j \alpha_{\sigma_{[r+n]}}\right\rvert
\\&+\sum_{j = J_0 + 1}^{\infty} \beta_{i+j}\left(
q_i \delta_{i+j} + 
C([r+j])\prod_{n=1}^j \alpha_{\sigma_{[r+n]}}\right).
\end{split}
\end{equation}
By \eqref{exp_grow}, we have $\frac{q_i}{q_{i+j}} \ll c^j$
for some $c < 1$ and by \eqref{qk_delta_k}, $q_{i+j}\delta_{i+j} \ls 1$. 
Similarly, we have $C([j+i])\prod_{n=1}^j \alpha_{\sigma_{[r+n]}}\ll \tilde{c}^j$
holds for some $0 < \tilde{c} < 1$, and since
$\alpha$ is badly approximable, $b_{i+j} \ll 1$,
where all implied constants only depend on $\alpha$.
Thus, we can follow that

\[\sum_{j = J_0 + 1}^{\infty} \beta_{i+j}\left(
q_i \delta_{i+j} + 
C([r+j])\prod_{n=1}^j \alpha_{\sigma_{[r+n]}}\right)
\ll \sum_{j = J_0+1}^{\infty} \left(c^j + \tilde{c}^j\right).
\]
Choosing $J_0(\delta)$ sufficiently large, we see that the second sum in \eqref{I0J0} is bounded from above by $\delta/2$.
Applying \eqref{quadr_limit1},
we have for every $j \gs 0$ that
\begin{equation*}
\begin{split} 
\lim_{m \to \infty} (-1)^{i}q_{m\ell + i} \delta_{m\ell + i+j}
= (-1)^j C([j+i])\prod_{n=1}^j \alpha_{\sigma_{[r+n]}}.
\end{split}
\end{equation*}
So if $i \gs I_0$ with $I_0 = I_0(J_0)$ sufficiently large, 
also the first sum in \eqref{I0J0} is bounded from above by $\delta/2$, which concludes the proof.
\end{proof}

\begin{proof}[Proof of Theorem \ref{l=2}]
By Theorem \ref{K_l<=7} and \eqref{liminf_limsup_equiv}, it suffices to show that $\liminf_{N \to \infty} P_N(\alpha) = 0$ holds for
those $\alpha$ where
$T^p(\alpha) = [0;\overline{6,a}]$ for some $p \in \mathbb{N}$ and $a \in \{1,2,3,4,5\}$. If $a \ls 4$,
then a direct computation of $G_0(T^p(\alpha),\cdot)$
shows $G_0(T^p(\alpha),0) < 1$ (see Figure \ref{6a})
and thus, the statement follows from Theorem B since this criterion is invariant under the application of the Gauss map. 
\captionsetup[subfigure]{labelformat=empty}
\begin{figure}[h]
\begin{subfigure}{.49\textwidth}
  \centering
  \includegraphics[width=.9\linewidth]{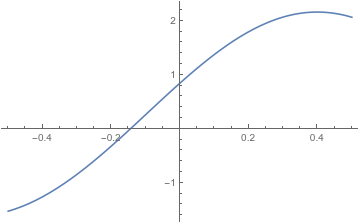}
  \caption{\phantom{ }\hspace{.7cm}
  $\alpha = [0;\overline{6,2}],\; G_0(\alpha,0) = 0.849\ldots$
  }
\end{subfigure}
\begin{subfigure}{.49\textwidth}
  \centering
  \includegraphics[width=.9\linewidth]{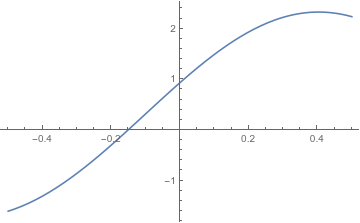}
  \caption{\phantom{ }\hspace{.7cm}
  $\alpha = [0;\overline{6,3}],\; G_0(\alpha,0) = 0.936\ldots$}
\end{subfigure}
\begin{subfigure}{.49\textwidth}
  \centering
  \includegraphics[width=.9\linewidth]{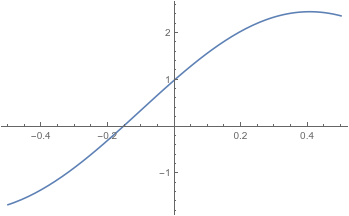}
  \caption{\phantom{ }\hspace{.7cm}
  $\alpha = [0;\overline{6,4}],\; G_0(\alpha,0) =0.998\ldots$
  }
\end{subfigure}
\begin{subfigure}{.49\textwidth}
  \centering
  \includegraphics[width=.9\linewidth]{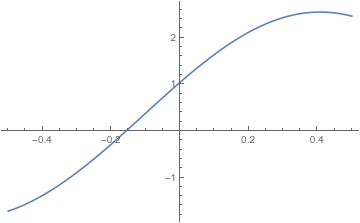}
  \caption{\phantom{ }\hspace{.7cm}
  $\alpha = [0;\overline{6,5}],\; G_0(\alpha,0) =1.047\ldots$}
\end{subfigure}
\caption{Limit function $G_0(\alpha,\varepsilon)$ for $\alpha = [0;\overline{6,a}]$ with $a \in \{1,2,3,4,5\}$. It is visible that for $a = 2$, we have $G_0(\alpha,0) < 1$. This inequality is also valid for $a = 3,4$, but the numerical value of $G_0(\alpha,0)$ gets closer to $1$ as $a$ increases. For $a = 5$, we obtain $G_0(\alpha,0) > 1$.
}
\label{6a}
\end{figure}
For the rest of the proof, we set 
$\alpha = [0;\overline{6,5}]$.
Now let
\[N_n = \sum_{i=0}^{2n} q_i(\alpha),\]
that is we consider only natural numbers with all Ostrowski coefficients being $1$.
Let $(\beta_j)_{j \in \mathbb{N}}$ with $\beta_j = 1$ for all $j \in \mathbb{N}$ (which is clearly an admissible sequence), and let $\delta > 0$
arbitrary.
By Proposition \ref{eps_convergence}, we have that for any $\delta$, there exist $I_0, J_0$ such that
\begin{equation*}
    \sup_{I_0 < i < 2n - J_0}  \left\lvert\varepsilon_{i,k}(N_n) - \varepsilon'_{[i],k}\left((\beta_{j})_{j \gs 1}\right) \right\rvert < \delta,
\end{equation*}
provided that $n$ is sufficiently large.
Using the notation $\alpha_{\pi} = {\alpha}_{\sigma_1}\cdot {\alpha}_{\sigma_2}$, we have
\[
\begin{split}
\varepsilon_{1,0}'\left((\beta_j)_{j \in \mathbb{N}}\right) &= - C(0)\frac{{\alpha}_{\sigma_2}}{1 - \alpha_{\pi}} + C(1)\frac{\alpha_{\pi}}{1 - \alpha_{\pi}}  \approx -0.025499,\\
\varepsilon_{0,0}'\left((\beta_j)_{j \in \mathbb{N}}\right) &= - C(1)\frac{{\alpha}_{\sigma_1}}{1 - \alpha_{\pi}} + C(0)\frac{\alpha_{\pi}}{1 - \alpha_{\pi}} \approx -0.0266289,
\end{split}
\]
which follows from elementary geometric series identities.
Using
\eqref{conv_of_eps}, we can therefore prove with computational assistance that
  \begin{equation*}
 G_0\left(\varepsilon_{2i,0}(N_n)\right)\cdot G_1\left(\varepsilon_{2i-1,0}(N_n)\right) < 0.997,\end{equation*}
 provided that $J_0 \ls 2i \ls 2n - I_0$. By Theorem A, we see that for
 $i \gs I_0(\delta)$,
 \begin{equation*}
 \left\lvert P_{q_i}\left(\alpha,\varepsilon_{i,0}(N_n)\right) - G_{[i]}\left(\alpha,\varepsilon_{i,0}(N_n)\right)\right\rvert < \delta.\end{equation*}
Choosing $\delta$ sufficiently small, we thus obtain for $I_0(\delta) < 2i < 2n- J_0(\delta)$ 
(without loss of generality, $I_0,J_0$ are even)
that
\[P_{q_{2i}}\left(\alpha,\varepsilon_{2i,0}(N_n)\right)\cdot P_{q_{2i-1}}\left(\alpha,\varepsilon_{2i+1,0}(N_n)\right)
< 0.999.\]
Using the decomposition \eqref{shifted_prod_1} from Proposition \ref{prop_shifted}, this implies
\begin{equation*}
\begin{split}
    P_N(\alpha) < &\prod_{i =1}^{I_0/2} P_{q_{2i}}\left(\alpha,\varepsilon_{i,0}(N_n)\right)P_{q_{2i-1}}\left(\alpha,\varepsilon_{i+1,0}(N_n)\right)\times
     \\&\prod_{i = n - J_0/2}^{n} P_{q_{2i}}\left(\alpha,\varepsilon_{i,0}(N_n)\right)P_{q_{2i-1}}\left(\alpha,\varepsilon_{i+1,0}(N_n)\right)\cdot
    0.999^{n-I_0/2-J_0/2}.
    \end{split}
\end{equation*}
Since $I_0,J_0$ are finite and by Lemma \ref{gg1}, all shifted products are uniformly bounded.
letting $n \to \infty$ shows that $\liminf_{N \to \infty} P_N(\alpha) = 0$.
We are left to show that this argumentation also applies for any $\beta \in T^{-p}(\alpha)$, $p \in \mathbb{N}$.
By Lemma \ref{aist_borda_limit_fct}, we have that
$\lim_{i \to \infty} P_{q_{i}}(T^p(\alpha),\varepsilon) - P_{q_{i+p}}(\alpha,\varepsilon) = 0$
and a variation of Proposition \ref{eps_convergence} with indices shifted by $p$ is still valid, so
choosing $(\beta_j)_{j \in \mathbb{N}} = (1,1,1,\ldots)$ shows $\liminf_{N \to \infty} P_N(\beta) = 0$
by the same argumentation as above.
\end{proof}

\section{Proof of Theorem \ref{counterex_l3}}\label{counterex_sec}

The proof of Theorem \ref{counterex_l3} follows some ideas of \cite[Lemma 3]{tech_zaf}, although we need to incorporate the additional difficulty of period lengths $\ell > 1$. We start with the easier case where
$\alpha = [0;\overline{5,4}]$ first, and prove the statement for $\alpha = [0;\overline{6,5,5}]$
later with similar ideas, only stating the needed refinements.
We say that a function $\Pi_r: \mathbb{N}^j \to \mathbb{R}, 0 \ls r \ls \ell$ is a \textit{lower-approximation function} for $\alpha$ and $r$ if for any admissible sequence $(\beta_{j})_{j \in \mathbb{N}}$ and any $i \in \mathbb{N}$ that
fulfills $b_{i-1} = \beta_1, b_i = \beta_2, b_{i+j-1} = \beta_j, [i] = r$, 
we have

\begin{equation*}
\Pi_r(\beta_1,\ldots,\beta_j) \ls
\prod_{k=1}^{b-1}G_{r}\left(\varepsilon'_{r,k}((b_{i+1+j})_{j \in \mathbb{N}})\right)\cdot
\left(G_{[r-1]}\left(\varepsilon'_{r-1,k}((b_{i+j})_{j \in \mathbb{N}}))\right)\right)^{\mathds{1}_{[b_{i-1} \neq 0]}},
\end{equation*}
where $\varepsilon'_{r,k}((\beta_{j})_{j \in \mathbb{N}})$ is defined as in \eqref{def_eps'_seq}.
The following lemma shows us how we prove the theorem under the assumption of having good lower-approximation functions $\Pi_{[i]}$ whose evaluation exceeds $1$ for
most $i$. These functions will be explicitly constructed in Lemmas \ref{compact54} respectively \ref{numerics_compact} by estimates on the perturbation value.

\subsection{Period length $\ell =2$}

\begin{lem}\label{suffic_admi}
Let $\alpha$ be a quadratic irrational with period $\ell = 2$ and $\Pi_0,\Pi_1: \mathbb{N}^4 \to \mathbb{R}$ be lower-approximation functions for $\alpha$ that fulfill
\begin{equation}\label{Pi>1.01}\Pi_1(a,b,c,d) \cdot \Pi_0(b,c,d,e) >1.01\end{equation}
for any admissible $(a,b,c,d,e)$ (with respect to $\alpha$ and $0$) where $(a,b,c) \neq (0,0,0)$.
Then we have 
\[\liminf_{N \to \infty} P_N(\alpha) > 0, \quad \limsup_{N \to \infty} \frac{P_N(\alpha)}{N} < \infty.\]
\end{lem}

\begin{proof}
Recall from Proposition \ref{prop_shifted} that we can write
\begin{equation}\label{shifted_prod_3}
    P_N(\alpha) = P_{q_n}(\alpha) \cdot \prod_{i = 1}^n K_i(N)
\end{equation}
with

\begin{equation*}K_i(N) = \prod_{i=0}^{n}\prod_{c_i= 1}^{b_i-1} P_{q_i}\left(\alpha,\varepsilon_{i,c_i}(N)\right)
\cdot P_{q_{i-1}}\left(\alpha,\varepsilon_{i-1,0}(N)\right)^{\mathds{1}_{[b_{i-1} \neq 0]}}.\end{equation*}
Let
$\delta >0$ be a fixed small constant. 
Using the uniform continuity of $P_{q_i}$ from Theorem A and the assumption that $\Pi_i$ is a lower approximation function, we have
\begin{equation}\label{KI_epsclose}K_i(N) \gs \Pi_{[i]}\left(b_{i-1}(N),b_i(N),b_{i+1}(N),b_{i+2}(N)\right) - \delta\end{equation}
for all $i \gs I_0 = I_0(\delta)$.
If $b_{i-1} = b_i = 0$, all considered products are empty, hence

\begin{equation}
    \label{empty_prod_case}
    K_i(N) = \Pi_{[i]}\left(b_{i-1}(N),b_i(N),b_{i+1}(N),b_{i+2}(N)\right) = 1.
\end{equation}
Now let $[i] = 1$ and set
$a = b_{i-1}(N), b = b_i(N), c = b_{i+1}(N), d = b_{i+2}(N)$. Combining \eqref{KI_epsclose} and \eqref{empty_prod_case}, we have

\begin{equation*}
K_{i}(N)\cdot K_{i+1}(N)
\gs \left(\Pi_{[i]}(a,b,c,d) - \mathds{1}_{[(a,b) \neq (0,0)]}\delta\right)\cdot \left(\Pi_{[i+1]}(b,c,d,e) - \mathds{1}_{[(b,c)\neq (0,0)]}\delta\right).
\end{equation*}
If $(a,b,c) \neq (0,0,0)$, then combining \eqref{Pi>1.01} with the previous inequality
yields
\begin{equation}\label{K_I>1}K_{i}(N)\cdot K_{i+1}(N) \gs 1,\end{equation}
provided that $i \gs I_0$, and $\delta$ is chosen sufficiently small.
Clearly, \eqref{empty_prod_case} shows that the previous inequality also holds in case $(a,b,c) = (0,0,0)$.
Combining \eqref{shifted_prod_3} and \eqref{K_I>1}, we obtain
\[P_{N}(\alpha) \gs P_{q_n}(\alpha)\cdot \prod_{i = 1}^{I_0} K_{i}(N).\]
By Lemma \ref{gg1}, each of those finitely many factors is uniformly bounded away from $0$, so clearly, we have $\liminf_{N \to \infty} P_N(\alpha) > C(I_0,\alpha) > 0$
 and by \eqref{liminf_limsup_equiv},
the result follows.
\end{proof}

So we are left to show that the assumption of Lemma \ref{suffic_admi} is fulfilled
when choosing $\alpha = [0;\overline{5,4}]$, which will be shown in the following statements.

\begin{lem}\label{compact54}
Let $\alpha = [0;\overline{5,4}]$, $(a,b) \neq (0,0)$ and let $(a,b,c,d)$ be admissible values with respect to $0$ respectively $1$. Then there exist lower-approximation functions $\Pi_r$
with respect to $\alpha$ and $r = 0,1$ with the following properties:
 \begin{enumerate}[(i)]
 \item $\Pi_0(a,b,c,d) > 1.01$ for all $(a,b) \neq (0,0)$. If $a \neq 0$, then $\Pi_0(a,b,c,d) > 1.22$.
     \item If $a = 0, b \neq 0$, then $\Pi_1(a,b,c,d) > 1.01$.
     \item If $a \neq 0$ and $b \neq 1$, we have $\Pi_1(a,b,c,d) > 1.01$.
     \item If $a \neq 0, b =1$, we have $\Pi_1(a,b,c,d) > 0.84$.
     \item If $(a,b) = (0,0)$, then $\Pi_0(a,b,c,d)= \Pi_1(a,b,c,d) = 1$.
 \end{enumerate}
\end{lem}
Assuming for the moment Lemma \ref{compact54} to hold, we can prove 
the following corollary.

\begin{cor}\label{prove_54}
Let $\alpha = [0;\overline{5,4}]$. Then we have

\[\liminf_{N \to \infty} P_N(\alpha) > 0,\quad 
\limsup_{N \to \infty} \frac{P_N(\alpha)}{N} < \infty.\]

\end{cor}

\begin{proof}
In view of Lemma \ref{suffic_admi}, it suffices to prove that for 
all admissible $(a,b,c,d,e)$ where $(a,b,c) \neq (0,0,0)$, we have
\begin{equation}\label{PI101}\Pi_1(a,b,c,d) \cdot \Pi_0(b,c,d,e) >1.01.\end{equation}
If $a = 0$, Lemma \ref{compact54}(ii) and (v) implies
$\Pi_1(a,b,c,d) \gs 1$. 
If $b \neq 0$ we have by Lemma \ref{compact54}(i) that
$\Pi_0(b,c,d,e) > 1.22$, so \eqref{PI101} clearly holds. If $(a,b) = (0,0)$, then $c \neq 0$,
hence by Lemma \ref{compact54}(i) and (v), we have \eqref{PI101} again.
Since $\Pi_0(b,c,d,e) \gs 1$ in any case, Lemma \ref{compact54}(iii) treats all cases except $b = 1$. But if $b = 1$, Lemma \ref{compact54}(i) and (iv) shows that
$\Pi_1(a,b,c,d)\cdot \Pi_0(b,c,d,e) \gs 1.22\cdot 0.84 > 1$.
\end{proof}

\begin{proof}[Proof of Lemma \ref{compact54}]
The main task in this proof is to construct good lower-approximation functions $\Pi_r$
and check with computational assistance that the stated numerical inequalities are fulfilled.
Similarly to Proposition \ref{eps_convergence}, we define for $r \in \{0,1\}$ and $1 \ls k
\ls a_{r+1}$
\[\varepsilon_{r,k}'(a,b)
:= kC(r) - aC([r+1]){\alpha}_{\sigma_[r+1]} + bC([r+2])\alpha_{\pi}.\]
For shorter notation, we write 
\[c_{r,0} = C(r), \quad c_{r,1} = C([r+1]){\alpha}_{\sigma_[r+1]}, \quad c_{r,2} = C([r])\alpha_{\pi}.
\]
Let $\varepsilon'_{r,k}((b_{i+j})_{j \in \mathbb{N}})$ be defined as in Proposition \ref{eps_convergence}, $[i] = r$ with $b_{i+1} = a, b_{i+2} = b$ and $(b_{i+j})_{j \in \mathbb{N}}$ a sequence of admissible values.
Observe that

\begin{equation*}
\begin{split}
 \varepsilon'_{r,k}((b_{i+j})_{j \in \mathbb{N}}) - \varepsilon_{r,k}'(a,b)
 = \; &c_{r,1}\left( b_{i+3}\alpha_{\pi} + b_{i+5}\alpha_{\pi}^2 + b_{i+7}\alpha_{\pi}^3 + \ldots\right)
 \\\;+ &c_{r,2}\left(b_{i+4}\alpha_{\pi} + b_{i+6}\alpha_{\pi}^2 + b_{i+8}\alpha_{\pi}^3 + \ldots\right).
\end{split}
\end{equation*}
Since $0 \ls b_{i+j} \ls a_{i+j+1} \ls 5$, we get
\begin{equation}\label{bounds_eps54}L_r := - 5c_{r,1} \frac{\alpha_{\pi}}{1 - \alpha_{\pi}} \ls \varepsilon'_{r,k}((b_{i+j})_{j \in \mathbb{N}}) - \varepsilon_{r,k}'(a,b) \ls 5c_{r,2} \frac{\alpha_{\pi}}{1 - \alpha_{\pi}} := U_r.\end{equation}
Now we define

\[
\tilde{G}_{r,k}(a,b) := \min\left\{G_{r}\left(\varepsilon'_{r,k}(a,b) + L_r\right),
G_{r}\left(\varepsilon'_{r,k}(a,b) + U_r\right)\right\}
\]
and

\begin{equation*}
\Pi_r(a,b,c,d) := \prod_{k=1}^{b-1}\tilde{G}_{r,k}(c,d)\cdot
\left(\tilde{G}_{[r-1],0}(b,c)\right)^{\mathds{1}_{[a \neq 0]}}.\end{equation*}
By Lemma \ref{approx_limit}(i), we have

\[
\Pi_r(a,b,c,d) \gs 
\prod_{k=1}^{b-1}G_{r}( \varepsilon'_{r,k}((b_{i+j})_{j \in \mathbb{N}}))\cdot
\left(G_{[r-1]}( \varepsilon'_{r,k}((b_{i+j-1})_{j \in \mathbb{N}}))\right)^{\mathds{1}_{[b_{i+j-1} \neq 0]}}
\]
for any admissible sequence $(\beta_{j})_{j \in \mathbb{N}}$ that fulfills
$\beta_{i-1} = a, \beta_i = b, \beta_{i+1} = c$ and $\beta_{i+2} = d$, hence $\Pi_0, \Pi_1$ are lower-approximation functions. 
Applying the procedure described in Remark \ref{explicit_com},
we can compute lower bounds for $\Pi_r(a,b,c,d)$ in finitely many steps. In that way, we obtain the results (i) - (v) with computational assistance by distinguishing all admissible cases for $a,b,c,d,e$.
This concludes the proof of Lemma \ref{compact54}.
\end{proof}

\subsection{Period length $\ell =3$}
Having established the result for $\alpha = [0;\overline{5,4}]$, we
head over to the case where $\alpha = [0;\overline{6,5,5}]$. The argument is very similar to the case above, so we will only indicate the changes that are needed.
We have the following statement, which is analogous to Lemma \ref{compact54}, followed by a corollary comparable to Corollary \ref{prove_54}:
\begin{lem}\label{numerics_compact}
Let $\alpha = [0;\overline{6,5,5}]$, and let $\Pi_r: \mathbb{N}^5 \to \mathbb{R}$ be lower-approximation functions for $\alpha$ and $r = 0,1,2$. Let $(a,b,c,d,e)$ be admissible values.
Then we have the following:
\begin{enumerate}[(i)]
    \item We have $\Pi_0(a,b,c,d,e)> 1.001, \Pi_1(a,b,c,d,e) > 0.81, \Pi_2(a,b,c,d,e) > 1.001$ for all admissible $(a,b,c,d,e)$ with $(a,b) \neq (0,0)$. If $b \neq 1$ and $a \neq 0$, then $\Pi_1(a,b,c,d,e) \gs 1.001$.
    \item If $b \neq 1$ and $a \neq 0, \Pi_2(a,b,c,d,e) \gs 1.25$, $\Pi_0(a,b,c,d,e) \gs 1.19$
    \item If $a \neq 0, c \gs 1$, then  $\Pi_1(a,b,c,d,e) \gs 0.85$.
    \item If $a \neq 0$ and $e \neq 0$ or $f \ls 2$, we have
    $\Pi_1(a,1,1,1,e)\cdot \Pi_2(1,1,1,e,f)\cdot \Pi_0(1,1,e,f,g) > 1.007.$
    \item If $a \neq 0$, we have 
    $\Pi_1(a,1,1,1,0)\cdot \Pi_2(1,1,1,0,f)\cdot \Pi_0(1,1,0,f,g) > 0.98.$
    \item If $(a,b) = (0,0)$, then 
    $\Pi_0(a,b,c,d,e) = \Pi_1(a,b,c,d,e) = \Pi_2(a,b,c,d,e) = 1$.
    \end{enumerate}
\end{lem}

\begin{cor}\label{cor_655}
\begin{enumerate}[(i)]
\item Let $a,b,c,d,e,f,g$ be admissible values such that $(a,b,c,d) \neq (0,0,0,0)$ and
\begin{equation}\label{prod<1}\Pi_1(a,b,c,d,e)\cdot \Pi_2(b,c,d,e,f)\cdot \Pi_0(c,d,e,f,g)< 1.001.\end{equation}
Then we have
$a \neq 0, b = 1, c = 1, d = 1, e = 0, f \gs 3$
and in that case,
\[\Pi_1(a,1,1,1,0)\cdot \Pi_2(1,1,1,0,f)\cdot \Pi_0(1,1,0,f,g)> 0.98.\]
\item
Let $(a,b,c,d,e,f,g,h,i,j)$ be admissible values
such that
\[\Pi_1(a,b,c,d,e)\cdot \Pi_2(b,c,d,e,f)\cdot \Pi_0(c,d,e,f,g)< 1.\]
Then we have
\[\Pi_1(a,b,c,d,e)\cdot \Pi_2(b,c,d,e,f)\cdot \Pi_0(c,d,e,f,g)\cdot \Pi_1(d,e,f,g,h)\cdot \Pi_2(e,f,g,h,i) \cdot \Pi_0(f,g,h,i,j)> 1.\]
\end{enumerate}
\end{cor}
\begin{proof}
Let $a,b,c,d,e,f,g$ be admissible values such that \eqref{prod<1} holds.
We see from Lemma \ref{numerics_compact}(i) and (vi) immediately that $a \neq 0, b = 1$. If $c \neq 1$, then by (ii) $\Pi_2(b,c,d,e,f) > 1.25$,
hence \[\Pi_1(a,b,c,d,e)\cdot \Pi_2(b,c,d,e,f)\cdot \Pi_0(c,d,e,f,g)
> 0.81\cdot 1.25 > 1.\] Thus $c = 1$
follows, which implies by (iii) that $\Pi_1(a,b,c,d) \gs 0.85$. If $d \neq 1$, we have by (ii) that \[\Pi_1(a,b,c,d,e)\cdot \Pi_2(b,c,d,e,f)\cdot \Pi_0(c,d,e,f,g)\gs 0.85 \cdot 1.19>1,\] hence $d = 1$. By (iv), this implies $e = 0$ and $f \gs 3$. Using (v), this concludes the first statement. For the second statement, we use part (i) of the corollary to see that $ e= 0, f \gs 3$. This implies by $(i)$ that $P_1(d,e,f,g) \gs 1, P_0(f,g,h,i) \gs 1$ and by (ii), $P_2(e,f,g,h) \gs 1.25$. Since $0.98\cdot 1.25 > 1$, the result follows.
\end{proof}

\begin{lem}\label{prove_655}
Let $\alpha = [0;\overline{6,5,5}]$. Then we have
\[\liminf_{N \to \infty} P_N(\alpha) > 0,\quad 
\limsup_{N \to \infty} \frac{P_N(\alpha)}{N} < \infty.\]
\end{lem}

\begin{proof}
Let $N = \sum_{i=0}^{n} b_{i}q_i(\alpha)$ be the Ostrowski expansion of some arbitrary integer $N$.
Now let
$M_i(N) := K_{3i+1}(N)\cdot K_{3i+2}(N)\cdot K_{3i+3}(N)$ with $K_i$ defined as in \eqref{def_Ki}.
Arguing as in Lemma \ref{suffic_admi}, it suffices to show that for $i \gs I_0$ we have
$\prod\limits_{i = I_0}^{m}
M_i(N) \gs C > 0
$
where
$m := \left\lfloor \frac{n}{3}\right\rfloor-1$ and $C = C(\alpha)$ is independent of $N$.
We define
\[J := \left\{j \in \{I_0+1,\ldots,m\}: M_J(N) < 1\right\}, \quad 
J-1 := \{j -1: j \in J\}.
\]
Choosing $I_0$ sufficiently large, we can apply Corollary \ref{cor_655} to deduce
$M_j(N)\cdot M_{j+1}(N) > 1$ for all $j \in J$, $M_i(N) \gs 1$ for all $i \gs I_0+1$
and $M_{I_0}(N) > 0.98.$
Hence we obtain
\[\prod_{i \gs I_0}^{m}
M_i(N) \gs 0.99 \cdot \prod_{j \in J}\left(M_j(N)\cdot M_{j+1}(N)\right) \cdot \prod_{\substack{i = I_0\\ i, i+1 \notin J}}^{m}M_i(N) \gs 0.98 
\]
and thus, the statement follows.
\end{proof}

\begin{proof}[Proof of Lemma \ref{numerics_compact}]
As in Lemma \ref{compact54}, we define for $r \in \{0,1\}$ and $1 \ls k
\ls a_{r+1}$
\[\varepsilon_{r,k}'(a,b,c)
:= kc_{r,0} - ac_{r,1} + bc_{r,2} - cc_{r,3},
\]
where

\[\begin{split}
&\alpha_{\pi} := {\alpha}_{\sigma{[r+1]}}{\alpha}_{\sigma{[r+2]}}{\alpha}_{\sigma{[r+3]}}, \quad
c_{r,0} = C(r), \quad c_{r,1} = C([r+1]){\alpha}_{\sigma{[r+1]}}, \\&c_{r,2} = C([r+2]){\alpha}_{\sigma{[r+1]}}{\alpha}_{\sigma{[r+2]}}, \quad 
c_{r,3} = C(r)\alpha_{\pi}.\end{split}
\]
Analogously to \eqref{bounds_eps54}, we can prove that 

\begin{equation*}
\begin{split}L_r := -\frac{6\alpha_{\pi}}{1 - \alpha_{\pi}^2}\left(\alpha_{\pi}c_{j,1} + c_{j,2} + \alpha_{\pi}c_{j,3}\right) &\ls \varepsilon'_{r,k}((b_{i+j})_{j \in \mathbb{N}}) - \varepsilon_{r,k}'(a,b,c) \\&\ls \frac{6\alpha_{\pi}}{1 - \alpha_{\pi}^2}\left(c_{j,1} + \alpha_{\pi}c_{j,2} + c_{j,3}\right) := U_r\end{split}\end{equation*}
and define

\[\begin{split}
\tilde{G}_{r,k}(a,b,c) &:= \min\left\{G_{r}\left(\varepsilon'_{r,k}(a,b,c) + L_r\right),
G_{r}\left(\varepsilon'_{r,k}(a,b,c) + U_r\right)\right\},
\\\Pi_r(a,b,c,d,e) &:= \prod_{k=1}^{b-1}\tilde{G}_{r,k}(c,d,e)\cdot
\left(\tilde{G}_{[r-1],0}(b,c,d)\right)^{\mathds{1}_{[a \neq 0]}}.\end{split}
\]
As before, we can deduce that $\Pi_r, r \in \{0,1,2\}$ are lower-approximating functions, and we can finish the proof as in Lemma \ref{compact54} by using computational assistance.
\end{proof}

\subsection*{Acknowledgements}
The author is grateful to Christoph Aistleitner for various comments on an earlier version of this paper.

\end{document}